\documentclass[11pt,a4paper]{article}

\usepackage[a4paper]{geometry}
\usepackage[T1]{fontenc}
\usepackage{amsfonts,amsthm,amssymb,amsmath,mathdots, bbm,mathabx,mathrsfs}

\usepackage{url}
\usepackage[latin1]{inputenc}
\usepackage{pdfpages}
\usepackage{epstopdf}
\usepackage{graphicx}
\usepackage{epsfig}
\usepackage{subfigure}
\usepackage{cite}
\usepackage[colorlinks]{hyperref}
\hypersetup{
linkcolor=blue,
citecolor=blue,
}
\usepackage{dsfont}

\newcommand{\E}{{\mathbb E}}
\newcommand{\R}{{\mathbb R}}

\renewcommand{\P}{{\mathbb P}}

\renewcommand{\L}{\mathbb L}

\newcommand{\PP}{{\mathcal P}}
\DeclareMathOperator*{\argmin}{arg\,min}
\newcommand{\uargmin}[1]{\underset{#1}{\argmin}\;}

\def\var{\mathop{\rm Var}\nolimits}%

\newcommand{\thefont}[2]{\fontsize{#1}{#2}\fontshape{n}\selectfont}
\newcommand{\1}{\rlap{\thefont{10pt}{12pt}1}\kern.16em\rlap{\thefont{11pt}{13.2pt}1}\kern.4em}

\newtheorem{thm}{Theorem}[section]
\newtheorem{rem}{Remark}[section]

\newtheorem{defi}{Definition}[section]

\newcommand{\br}{\boldsymbol{r}}
\newcommand{\bq}{\boldsymbol{q}}
\newcommand{\bnu}{\boldsymbol{\nu}}
\newcommand{\bmu}{\boldsymbol{\mu}}
\newcommand{\bX}{\boldsymbol{X}}
\newcommand{\bF}{\boldsymbol{F}}
\newcommand{\bfun}{\boldsymbol{f}}

\title{Statistical data analysis in the Wasserstein space}
\author{J\'{e}r\'{e}mie Bigot \footnote{J\'er\'emie Bigot is a member of Institut Universitaire de France (IUF), and this work has been carried out with financial support from the IUF. I would like to thank Hervé Cardot and Pierre Calka for proposing to write this review, as well as the two anonymous referees for their comments and suggestions of improvement.} \\ Universit\'e de Bordeaux \\
Institut de Math\'ematiques de Bordeaux et CNRS  (UMR 5251)}

\begin{document}

\maketitle

\begin{abstract}
This paper is concerned by statistical inference problems from a data set whose elements may be modeled as random probability measures such as multiple histograms or point clouds. We propose to review recent contributions in statistics on the use of Wasserstein distances and tools from optimal transport to analyse such data. In particular, we highlight the benefits of using the notions of barycenter  and geodesic PCA in the Wasserstein space for  the purpose of learning  the principal modes of  geometric variation in a dataset. In this setting, we discuss existing works and we present some research perspectives related to the emerging field of statistical optimal transport.
\end{abstract}

\section{The emerging field of statistical optimal transport}

In many fields of interest (e.g.\ in signal and image processing or  bio-informatics), one records data in the form of high-dimensional vectors or matrices.  For the purpose of learning information form such data sets, a standard statistical analysis consists in considering that the observations are realizations of random variables belonging to a linear space endowed with an Euclidean distance. However, it has now been widely recognized that a significant gain in statistical inference can be achieved by the use of non-Euclidean distances to better capture the geometry of the sets to which the data (in the form of histograms, curves, images or points clouds) truly belong.  Such a choice of non-Euclidean distances comes from the observation that the data, in many applications, tend to vary along geodesics that are not straight lines. Therefore, the use of non-Euclidean distances is often meaningful to extract geometric information on non-linear sources of variability. Concepts such as nonlinearity and an appropriate modeling of the location of signals or images intensities are essential tools to solve high-dimensional regression problems for the purpose of estimation and classification. The use of Wasserstein distances based on the mathematics of optimal mass transport \cite{villani2003topics} allows to incorporate both spatial and intensity informations when comparing elements in a dataset that can be modeled as probability distributions. This approach leads to an intuitive geometric interpretation of the variability in such datasets.

\begin{figure}
\centering
\includegraphics[width= 0.6\textwidth,height=0.4\textwidth]{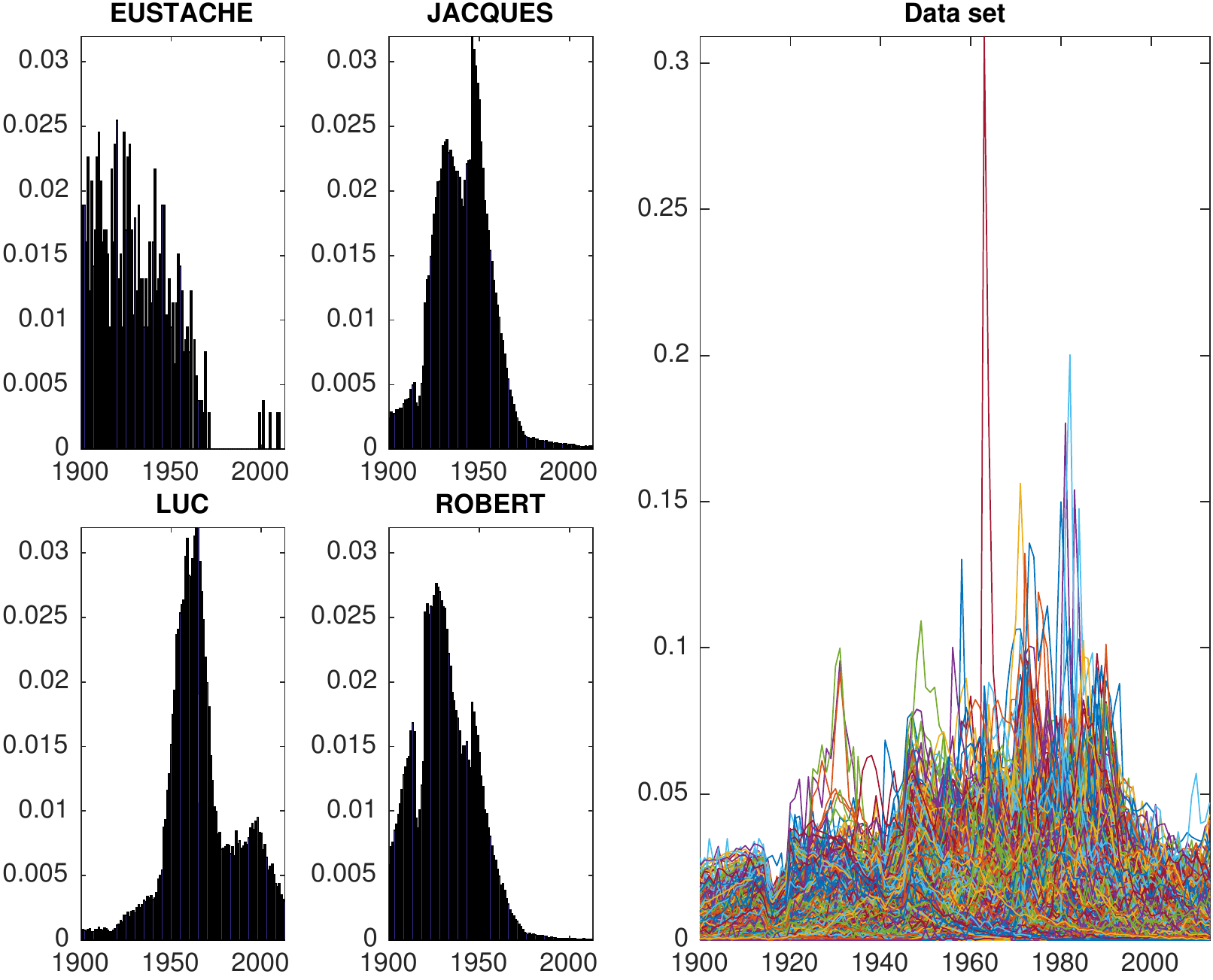} 
\caption{Children's first name at birth  - this dataset has been provided by the INSEE (French Institute of Statistics and Economic Studies). A subset of 4 histograms representing the distribution of children born with that name per year in France, and the whole dataset of $n=1060$ histograms (right), displayed as pdf over the interval $[1900,2013]$} \label{fig:names}
\end{figure}

\begin{figure}
\centering
\includegraphics[width=0.85 \textwidth,height=0.65\textwidth]{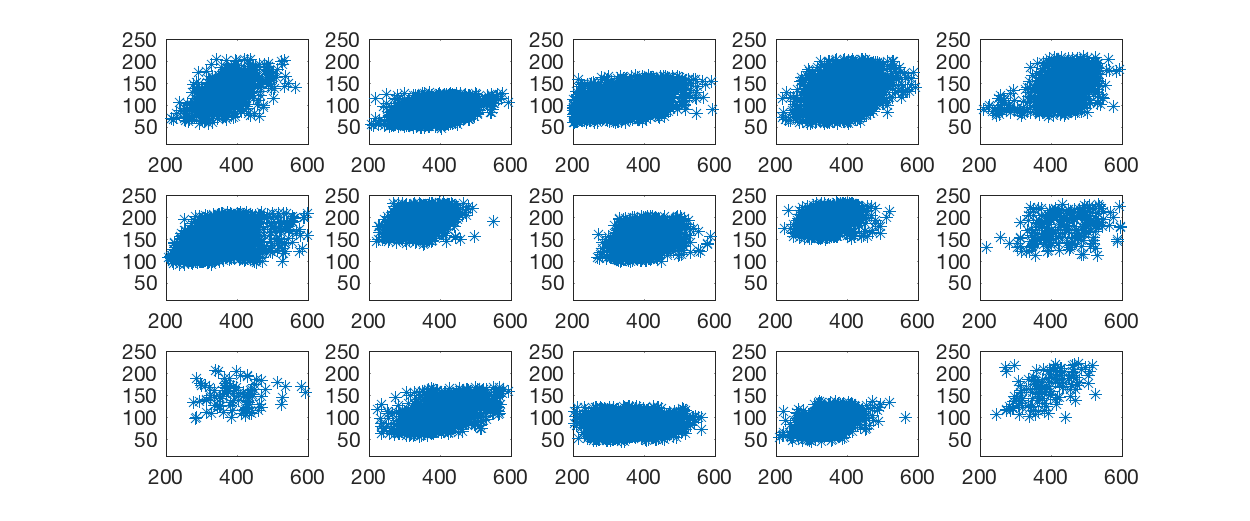} 
\caption{Example of  flow cytometry data \cite{hahne2010per} obtained from a renal transplant retrospective study conducted by the Immune Tolerance Network (ITN) and measured from $n=15$ patients (restricted to a bivariate projection).  The horizontal axis - resp.\  vertical axis represents the values of the forward-scattered light (FSC) -  resp.\ side-scattered light (SSC) cell marker, after an appropriate scaling trough  an arcsinh transformation and  an initial gating on total lymphocytes to remove artefacts. The number of considered cells by patient varies from $88$ to $2185$.}
\label{fig:ex_cytometry2D}
\end{figure}

Examples of such data are numerous. As an illustration, one may consider   histograms that represent, for a list of first names, the distribution of children born with that name per year in France between years 1900 and 2013. In Figure \ref{fig:names}, we display the histograms of four different names, as well as the whole  dataset. Another example arises from the statistical analysis of data acquired from flow cytometry which is a high-throughput technique in biotechnology that can measure a large number  of surface and intracellular markers of single cell in a biological sample. With this technique, one can assess individual characteristics (in the form of multivariate data) at a cellular level to determine the type of cell, their functionality and the way they differ. The development of this technology now leads to datasets made of multiple measurements (e.g.\ up to 18) of millions of individuals cells from different subjects. This leads to statistical inference problems from multiple points clouds such as those displayed in Figure \ref{fig:ex_cytometry2D} that can be modeled as discrete probability measures supported on $\R^d$.

Techniques based on optimal transport for data science have thus recently received an increasing interest in mathematical and computational statistics \cite{bigotbercu2019,bigot2018datadriven,BGKL18,BK17,BCP17,bigot2017geodesic,bigot2017central,DELBARRIO2019341,Pana15,Pana17,sommerfeld2016inference,spok2019,Klatt18,RIGOLLET20181228,NIPS2015_5680,delBarrio2019,alvarez-esteban2018,gouic2015existence,Petersen2019}, machine learning \cite{cuturi2013fast,rolet2016fast,NIPS2016_6566,pmlr-v84-genevay18a,cuturi,pmlr-v70-arjovsky17a,Frogner:2015,pmlr-v97-gordaliza19a,Flamary2018,Schmitz2018,NIPS2017_6792}, image processing and computer vision \cite{gramfort2015fast,rabin2015convex,benamou2015iterative,CuturiPeyre,ferradans2014regularized,2015-bonneel-siims,charlier,Solomon:2015,Solomon:2016} or computational biology \cite{article-cell}.

This paper, whose content is biased by the point of view of a statistician, is an introduction to some recent developments of optimal transport for statistical data analysis.  It is focused on the notions of barycenter and geodesic principal component analysis (PCA) in the Wasserstein space and their connection to nonparametric statistics and functional data analysis. We also report the results of numerical experiments to highlight, in an intuitive manner, the potential benefits of  such tools for statistical inference. A recent review on the statistical aspects of Wasserstein distances with a rich bibliography on the mathematics of optimal transport can also be found in \cite{Pana18}. More generally, computational optimal transport for data science is a rapidly growing research field. For a detailed introduction and examples of application beyond the prism of statistics we refer to \cite{bookOT}.

In Section \ref{sec:wass}, we present the notions of  Wasserstein space and  regularized optimal transport. Section \ref{sec:barW2} is focused on Wasserstein barycenters and  statistical inference problems related to this concept. Finally, geodesic principal component analysis (GPCA) in the Wasserstein space is briefly presented in  Section \ref{sec:GPCA}.  
Note that the presentation of the paper is rater informal on the mathematical aspects, and we mainly give intuition of the benefits of statistical optimal transport through the presentation of various numerical examples. Throughout the paper, we also discuss some research perspectives and open problems related to the statistical aspects of barycenters and GPCA in the Wasserstein space.

\section{Wasserstein space and regularized optimal transport} \label{sec:wass}

\subsection{The Wasserstein metric}

We denote by $\PP_2(\Omega)$  the space of probability measures with finite second moment and supports included in a convex domain $\Omega \subset \R^{d}$. This set can be endowed with the Wasserstein metric $W_2$ defined as
\begin{equation}
W_2(\mu,\nu)=\inf_{\pi\in\Gamma(\mu,\nu)}\left(\iint_{\Omega^2}\vert x-y\vert^2d\pi(x,y)\right)^{1/2}, \quad \mbox{ for } \mu,\nu \in \PP_2(\Omega), \label{def:W2}
\end{equation}
where $\Gamma(\mu,\nu)$ is the set of probability measures (also called transport plans or coupling measures) on the product space $\Omega\times\Omega$ with respective marginals $\mu$ and $\nu$, and $\vert \cdot \vert$ denotes the usual Euclidean norm on $\R^d$. The metric space $(\PP_2(\Omega), W_2)$ is called the ($2$-)Wasserstein space.

\begin{rem}
 The minimisation problem \eqref{def:W2} corresponds the so-called Kantorovich formulation of optimal transport. Other cost functions  than the quadratic one $c(x,y) = \vert x-y\vert^2$ may be used, but we restrict the discussion to this setting to simplify the presentation. For further details on the formulation of optimal transport for probability measures support on general metric space we refer to \cite{villani2003topics}.
\end{rem}

\begin{figure}
{\subfigure[Gaussian measure $\mu$ with density $f_{\mu}$]{\includegraphics[width=0.45\textwidth,height=0.25\textwidth]{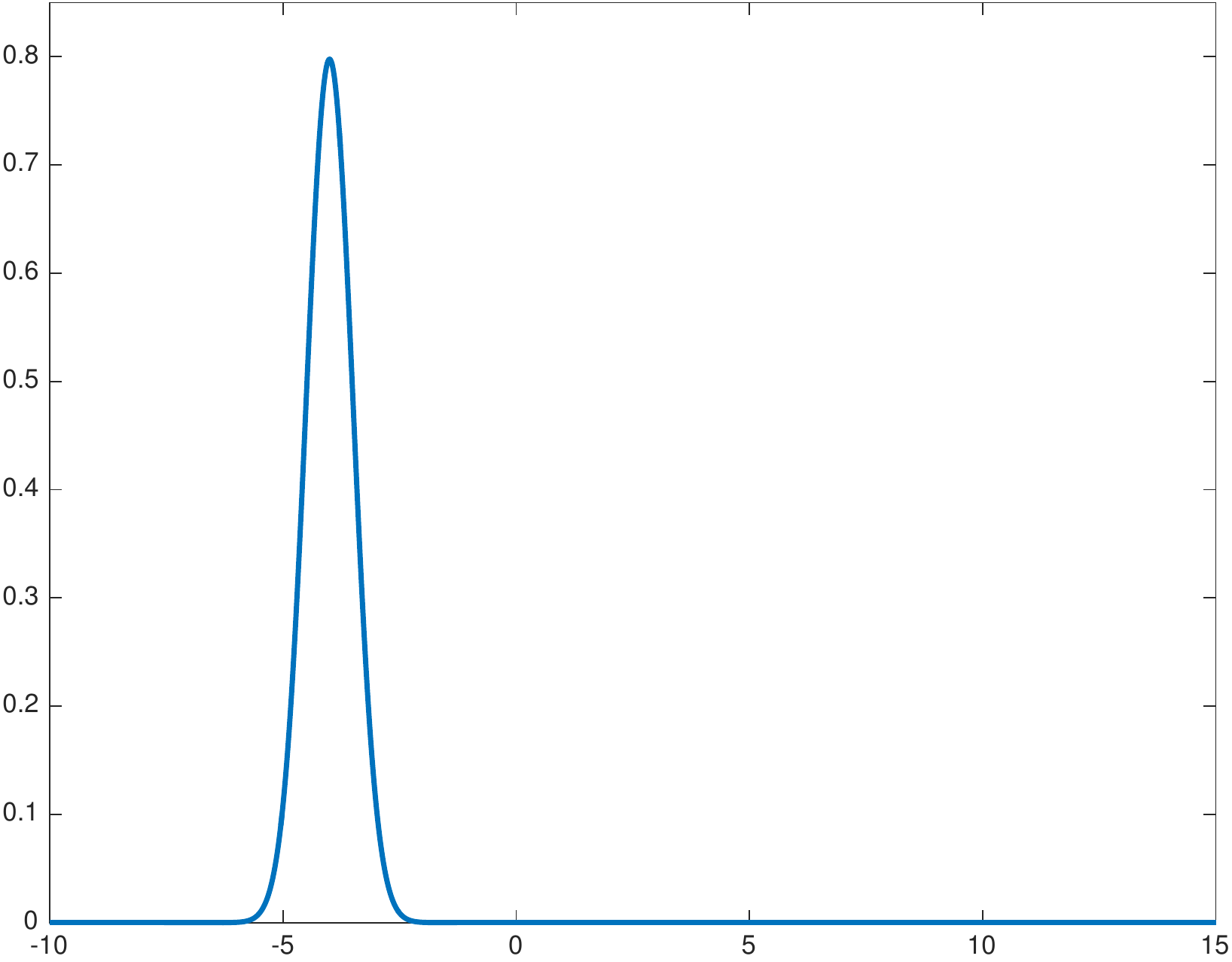}}}
{\subfigure[Gaussian measure $\nu$ with density $f_{\nu}$]{\includegraphics[width=0.45\textwidth,height=0.25\textwidth]{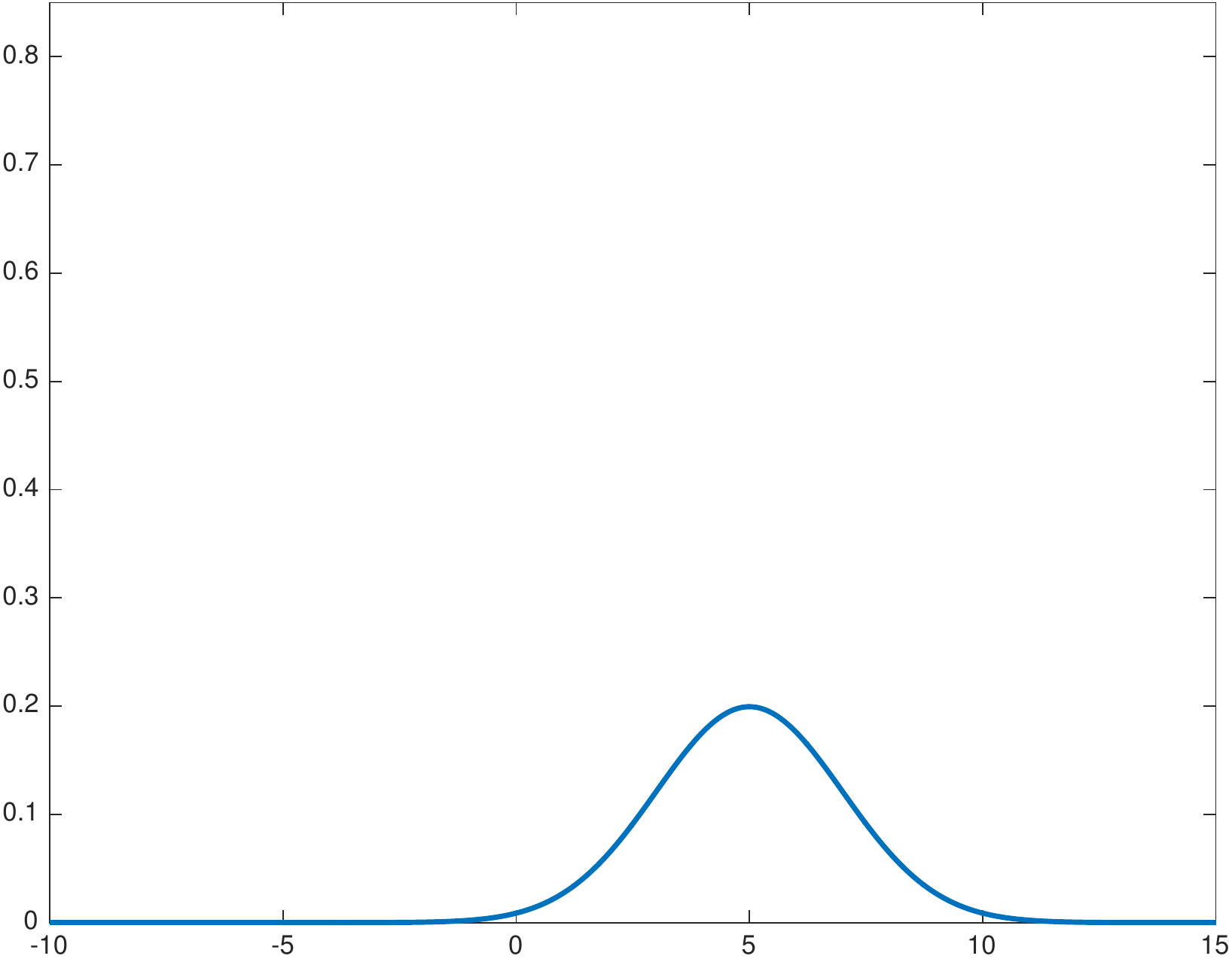}}}

{\subfigure[Geodesic in the Hilbert space of densities]{\includegraphics[width=0.45\textwidth,height=0.3\textwidth]{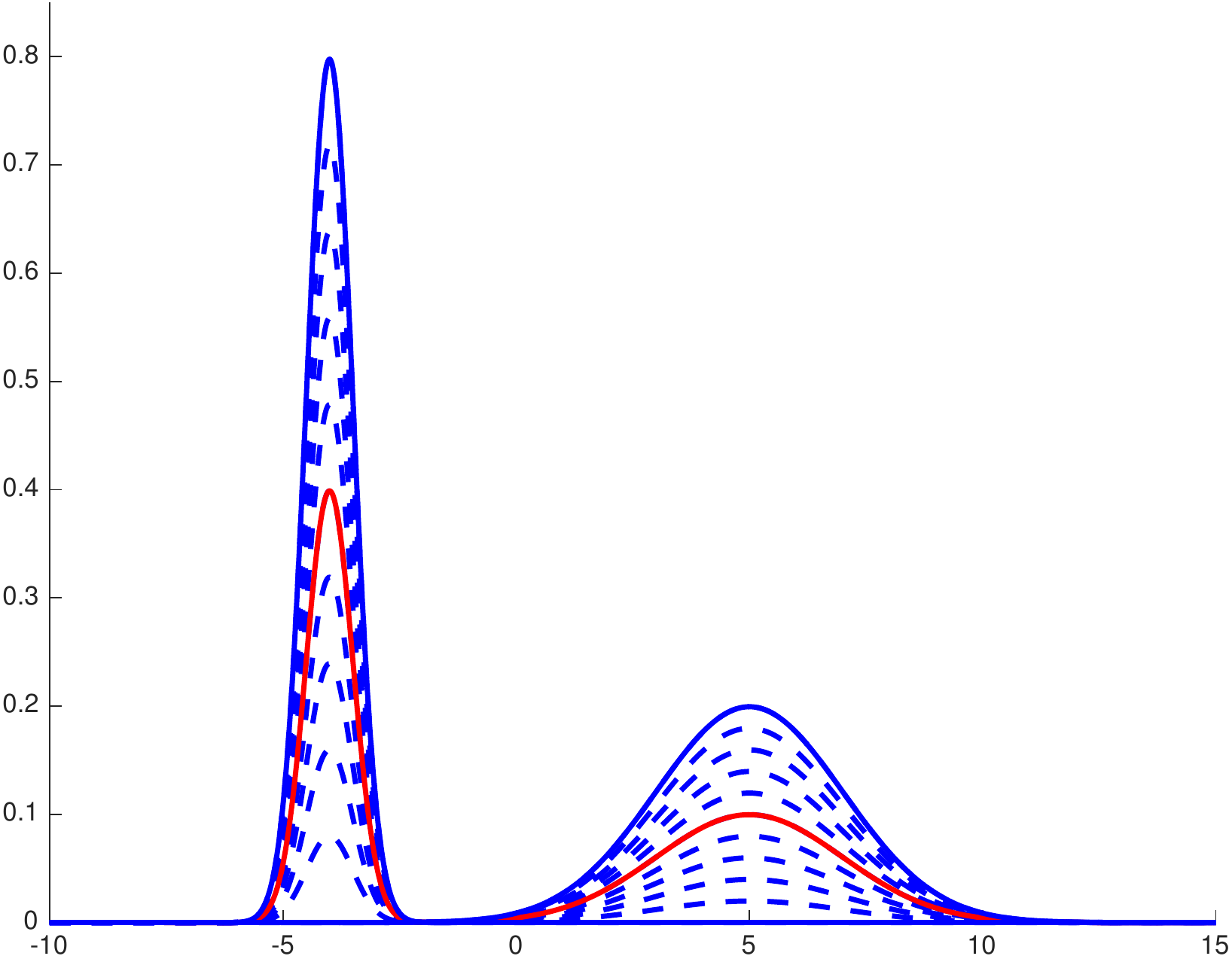}}}
{\subfigure[Geodesic in Wasserstein space]{\includegraphics[width=0.45\textwidth,height=0.3\textwidth]{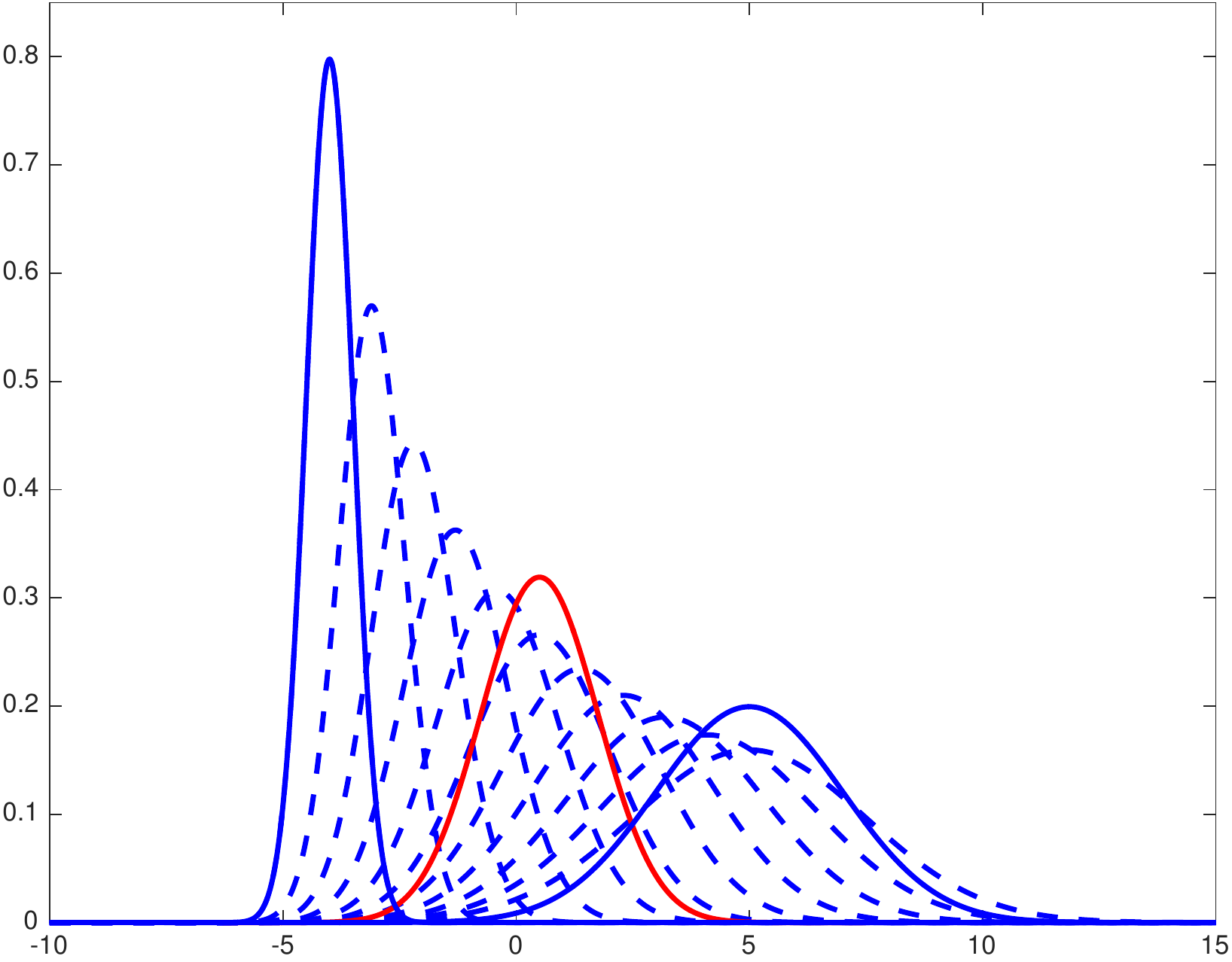}}}
\caption{First row: densities of two Gaussian measures on $\R$. Second row: geodesics between $\mu$ and $\nu$ with respect to the Hilbert metric and the Wasserstein metric. The dotted curves in (c) correspond to the densities $f_{t} = (1-t)f_{\mu} + t f_{\nu}$ , and in (d) they correspond to the densites of the measures with quantile functions $F^{-}_t = (1-t)F^{-}_{\mu} + tF^{-}_{\nu}$ for $0 < t < 1$. The red curves correspond to $t = 1/2$.} \label{fig:exintro1D}
\end{figure}

For probability measures supported on the real (that is when $d=1$), the Wasserstein metric has a closed form expression
\begin{equation}
W_2^2(\mu,\nu) = \int_{0}^{1} (F^{-}_{\mu}(u) - F^{-}_{\nu}(u))^{2} du,  \label{eq:W2_1d}
\end{equation}
for any $\mu,\nu \in \PP_2(\R)$,
where $F^{-}_{\mu}$ (resp.\ $F^{-}_{\nu}$) denotes the quantile function of the measure $\mu$ (resp.\ $\nu$). Hence, in the one-dimensional setting, the use of the Wasserstein metric amounts to perform data analysis in the space of quantile functions (supported on $[0,1]$) endowed with the usual Hilbert $\mathbb{L}_2$ metric. When $\mu$ and $\nu$ are absolutely continuous (a.c.) measures with square integrable probability density functions (pdf) $f_{\mu}$ and $f_{\nu}$, this remark allows to easily illustrate the difference between choosing either the Wasserstein metric $W_2(\mu,\nu)$ or the Hilbert metric $\|f_{\mu} - f_{\nu} \|_2 := \left( \int_0^1 |f_{\mu}(x) - f_{\nu}(x)|^2 dx\right)^{1/2}$ for  statistical analysis. Indeed, this choice implies two different constructions of the modes of variation in a dataset as illustrated in Figure \ref{fig:exintro1D}, where we display the geodesics between two Gaussian distributions in the Hilbert space of densities and in the Wasserstein space of probability measure. Recall that, in the one-dimensional setting, the value of the geodesic at time $t \in [0,1]$ between $\mu$ and $\nu$ (for the Wasserstein metric) is the measure $\mu_t$ whose quantile function is $F^{-}_{\mu_t} = (1-t) F^{-}_{\mu} + t F^{-}_{\nu}$. Therefore, the measure $\mu_t$  remains Gaussian along the geodesic path between two Gaussian measures $\mu$ and $\nu$. This illustrates the fact that performing statistical data analysis the Wasserstein space allows  to simultaneously take into account  variations in location and intensity (or mass) in a data set of densities (or histograms). A similar behavior holds in higher-dimensions ($d \geq 2$).

\subsection{Regularized optimal transport} 

A serious limitation for the use of the Wasserstein metric for data analysis is its computational cost between probability measures supported in higher-dimensional spaces (that is beyond the case $d=1$). Indeed,  the computational cost to evaluate a Wasserstein metric between two discrete probability distributions with supports of equal size $N$  is generally of order $\mathcal{O}(N^3\log N)$.   In this discrete setting, computing a Wasserstein distance amounts to solve the linear program  \eqref{def:W2} whose solution is constrained to belong to the convex set $\Gamma(\mu,\nu)$. Hence, the cost of this convex minimization becomes  prohibitive for moderate to large values of $N$. To reduce  the complexity of a linear program, a classical approach in optimization   \cite{wilson1969use} is to add a regularizing  entropy term. This is the approach followed in  \cite{cuturi,cuturi2013fast} by adding  entropy regularization to the transport transport plan in \eqref{def:W2}, which yields the strictly convex problem \eqref{def:primal_Sinkhorn}  below.

To simplify its presentation, we consider discrete probability measures supported on the same finite space $\Omega_{N}=\{x_1,\ldots,x_N\}\in(\R^d)^N$ of cardinal $N$. A discrete probability distribution $\mu$ (with fixed support included in $\Omega_{N}$) is identified by a vector of weights in the simplex $\Sigma_N=\{r=(r_1,\ldots,r_N)\in\R^N_{+} \ \mbox{with} \ \sum_{k=1}^Nr_k=1\}$ such that $\mu =\sum_{k=1}^Nr_k\delta_{x_k}$ where $\delta_x$ is the Dirac distribution at $x$. The Wasserstein distance between two discrete distributions $\mu =\sum_{k=1}^Nr_k\delta_{x_k}$ and $\nu =\sum_{k=1}^Nq_k\delta_{x_k}$ (identified to vectors $r$ and $q$ in $\Sigma_N$) then becomes
$$W_2(r,q):=\underset{U\in U(r,q)}{\min}\  \langle C, U\rangle^{1/2},$$
where  the set of coupling matrices is defined as $U(r,q):=\{U\in\R^{N\times N}_{+} \ \mbox{such that} \ U\mathds{1}_N=r,\ U^T\mathds{1}_N=q\}$ with $\mathds{1}_N$ the $N$ dimensional vector with all entries equal to $1$ and $C$ the cost matrix given by $C_{ml}=\vert x_m-x_l\vert^2$, for all $m,l\in\{1,\ldots,N\}$. Then, the notion of entropy regularized optimal transportation \cite{cuturi,CarlierDPS17} leads to the so-called Sinkhorn divergence $W_{2,\varepsilon}(r,q)$ between two discrete measures $r,q\in\Sigma_N$ which is defined as
\begin{equation}\label{def:primal_Sinkhorn}
W_{2,\varepsilon}^2(r,q):=\underset{U\in U(r,q)}{\min}\  \langle C, U\rangle -\varepsilon \textrm{h}(U).
\end{equation}
where  $\varepsilon>0$ is a regularization parameter, and the discrete (negative) entropy for a given coupling matrix $U\in U(r,q)$ is given by $ \textrm{h}(U):=-\sum_{m, \ell}U_{m \ell}\log U_{m \ell}$. The introduction of Sinkhorn divergence has been the starting point of the development of computational optimal transport \cite{bookOT} for machine learning as it makes feasible (from a computational point of view) the use of smoothed Wasserstein metrics for data analysis. Moreover, there now exists various toolbox and librairies to carry out data analysis using  (possibly regularized) optimal transport in Matlab - \url{https://optimaltransport.github.io/},  Python - \url{https://pot.readthedocs.io/en/stable/} and  the statistical computing environment R - \url{https://cran.r-project.org/package=Barycenter}.  

\section{Wasserstein barycenters} \label{sec:barW2}

\begin{figure}[htbp]
{\subfigure[$\nu_{1}$]{{\includegraphics[width=0.25 \textwidth,height=0.3\textwidth]{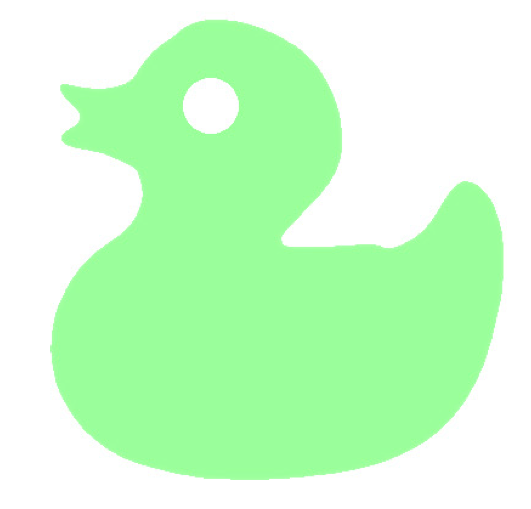}}}}
{\subfigure[$\bar{\mu}_2$]{{\includegraphics[width=0.25 \textwidth,height=0.3\textwidth]{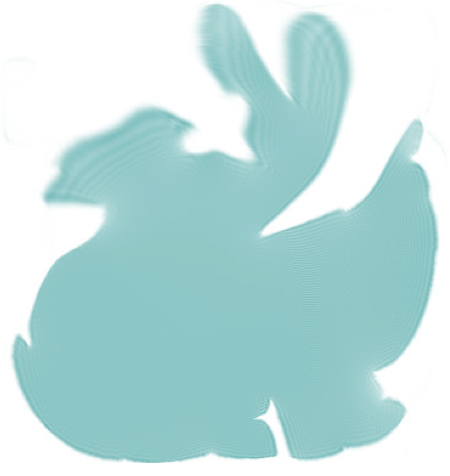}}}}
{\subfigure[$\nu_{2}$]{{\includegraphics[width=0.25 \textwidth,height=0.3\textwidth]{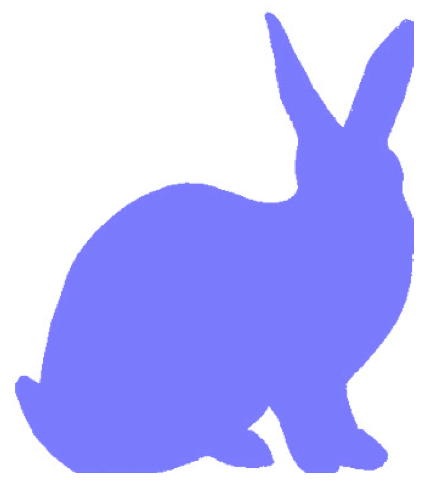}}}}
\caption{A numerical example taken from \cite{2015-bonneel-siims}. Wasserstein barycenter $\bar{\mu}_2$ in (b) between two probability measures $\nu_1$ and $\nu_2$ on $\R^2$ assumed to a.c.\ and uniform over their support displayed in (a) and (c). } \label{fig:exintro2D}
\end{figure}

The most basic statistical inference task from a data set is certainly to compute first order statistics that is to estimate an average location of the data. In this paper, we  focus on the estimation of  a population mean measure (or density function) from a dataset whose elements are probability measures. The notion of averaging depends on the metric that is chosen to compare elements in a given data set. In this paper, we consider the Wasserstein distance $W_2$ associated to the quadratic cost for the comparison of probability measures.  As introduced in \cite{agueh2011barycenters}, an empirical Wasserstein barycenter $\bar{\mu}_{n}$ of set of $n$ probability measures $\nu_{1},\ldots,\nu_{n}$ (not necessarily random) in $\PP_2(\Omega)$ is defined as a (possibly not unique) minimizer of
\begin{equation}
\mu \mapsto  \frac{1}{n} \sum_{i=1}^{n} W_2^2(\mu,\nu_{i}), \mbox{ over }  \mu \in \PP_2(\Omega). \label{eq:pbempBar}
\end{equation}
The Wasserstein barycenter  corresponds to the notion of empirical Fr\'echet mean \cite{fre} that is an extension of the usual Euclidean barycenter to nonlinear metric spaces. For the case $d=1$, the quantile formula \eqref{eq:W2_1d} implies that $\bar{\mu}_{n}$ is the measure with quantile function
$
F^{-}_{\bar{\mu}_n} = \frac{1}{n} \sum_{i = 1}^{n} F^{-}_{\nu_i}.
$
 For $n=2$, the Wasserstein barycenter  $\bar{\mu}_{2}$ is the probability measure located at half distance (that is at ``time $t=1/2$'') along the geodesic between two absolutely continuous (a.c.) measures $\nu_1$ and $\nu_2$. As illustrative examples, we display Wasserstein barycenters  $\bar{\mu}_{2}$ between two probability measures in dimension $d=1$  in Figure \ref{fig:exintro1D} and $d=2$ in Figure \ref{fig:exintro2D}.

Computing a Wasserstein barycenter between $n \geq 3$ probability measures is generally a delicate optimisation problem in dimension $d \geq 2$ that we do not discuss here. We refer to \cite{bookOT,benamou2015iterative,CuturiPeyre} for further details and references.  In what follows, we  focus on the mathematical statistical aspects of Wasserstein barycenters of random probability measures and the problem of their estimation from a nonparametric point of view. 

\subsection{Statistical models of random probability measures}

Let us first introduce the notion of a population Wasserstein barycenter, as well as a two-sampling model of random probability measures that is appropriate to consider the problem of its nonparametric estimation. To this end,  let $\P$  be a distribution on the space of probability measures  $(\PP_2(\Omega), \mathcal{B} \left( \PP_2(\Omega)\right)$, where $\mathcal{B} \left( \PP_2(\Omega) \right)$ is the  Borel $\sigma$-algebra  generated by the topology induced by  $W_{2}$. Throughout the paper, it is always assumed that $\P$ satisfies the square-integrability condition
$$
\int_{\PP_2(\Omega)}W_2^2(\mu,\nu)d\P(\nu) < + \infty, \mbox{ for some } \mu \in  \PP_2(\Omega).
$$
\begin{defi}
A {\bf population Wasserstein barycenter} $\mu_{\P}$ of a distribution $\P$ over $\PP_2(\Omega)$ is defined as  a (possibly not unique) minimizer of
\begin{equation}
\mu \mapsto  \int_{\PP_2(\Omega)}W_2^2(\mu,\nu)d\P(\nu), \mbox{ over }  \mu \in \PP_2(\Omega).  \label{eq:pbpopBar}
\end{equation}
\end{defi}

For a discrete distribution  $\P_{n} =\frac{1}{n}\sum\delta_{\nu_i}$ on $\PP_2(\Omega)$, one has that $\mu_{\P_{n}}$ corresponds to the empirical Wasserstein barycenter defined as a minimizer of \eqref{eq:pbempBar}. The existence and unicity of Wasserstein barycenters has been discussed in various works \cite{agueh2011barycenters,AGUEH2017812,gouic2015existence,BK17,KimPass}. As argued in \cite{AGUEH2017812},  a sufficient condition for the unicity of $\mu_{\P}$ is that the distribution $\P$ gives a strictly positive mass to the subset $\PP_2^{ac}(\Omega)$ of a.c.\ measures in $\PP_2(\Omega)$. Moreover, if $\P$ is supported on the set of measures $\PP_2(\Omega)\cap \mathbb{L}_q(\Omega)$ for some $q\in(1,+\infty)$ (i.e.  measures in $\PP_2^{ac}(\Omega)$ with $\mathbb{L}_q(\Omega)$ densities),  it follows that  $\mu_{\P}$  admits a density in $\mathbb{L}_q(\Omega)$. 
 
To define an appropriate sampling model, let us now assume that $\bnu_1,\ldots,\bnu_n$ are independent and identically distributed a.c.\ random measures  sampled from a  distribution $\P$ satisfying $\P\left( \PP_2^{ac}(\Omega) \cap \mathbb{L}_2(\Omega) \right) = 1$. Hence, under such assumptions, the population Wasserstein barycenter $\mu_{\P}$ is  unique and absolutely continuous. As raw data in the form of densities are generally not directly available, we focus on the setting where one has access to a set of  random vectors  $(\bX_{i,j})_{1 \leq j \leq p_{i}; \;  1 \leq i \leq n}$ in $\R^{d}$  organized in the form of $n$ experiment units, such that, conditionaly on $\bnu_i$, one has that   $\bX_{i,1},\ldots,\bX_{i,p_{i}}$ are iid observations sampled from the measure $\bnu_i$  for each $1 \leq i \leq n$. The observed data are thus in the form of multiple point clouds
$
\bX = (\bX_{i,j})_{1 \leq i \leq n; \; 1 \leq j \leq p_{i} }.
$
Elements in such a set can be modeled as $n$ discrete probability measures
\begin{equation}
\bnu_{p_{i}} = \frac{1}{p_{i}} \sum_{j=1}^{p_{i}} \delta_{\bX_{i,j}}, \label{eq:discrete}
\end{equation}
from which one may construct another empirical Wasserstein barycenter  $\hat{\bmu}_{n,p}$ defined as a (not necessarily unique)  minimizer of
\begin{equation}
\mu \mapsto  \frac{1}{n} \sum_{i=1}^{n} W_2^2(\mu,\bnu_{p_{i}}), \mbox{ over }  \mu \in \PP_2(\Omega), \label{eq:pbemp2Bar}
\end{equation}
with the convention that $p = \min_{1 \leq i \leq n} p_{i}$. We shall now discuss the issue of estimating  the population Wasserstein barycenter  $\mu_{\P}$ using either $\mu_{\P_{n}}$ (with $\P_{n} =\frac{1}{n}\sum\delta_{\bnu_i}$) or $\hat{\bmu}_{n,p}$ as $n \to +\infty$ and $p \to  +\infty$.

\subsection{Convergence of empirical Wasserstein barycenters}

\subsubsection{Law of large numbers}

By the law of large numbers, we refer to the question of showing that $W_2(\mu_{\P_{n}},\mu_{\P})$ almost surely converges to zero as $n \to + \infty$. This problem has been addressed in \cite{BK17} for a parametric model of random distributions $\P$ over $\PP_2(\Omega)$ with $\Omega$ assumed to be compact. General results for arbitrary distributions can be found in \cite{gouic2015existence}. Beyond the law of large numbers, it is natural to ask for a (functional) central limit theorem for the empirical Wasserstein barycenter. Nevertheless, this notion is more delicate to formulate as the Wasserstein space is not a Hilbert space. For recent contributions in this direction (but in specific parametric settings) we refer to \cite{AGUEH2017812,spok2019}.

\subsubsection{Rate of convergence in the one dimensional case}

Let us now restrict to the one dimensional case of probability measure in $\PP_2(\R)$ and assume that $p_1 = \ldots = p_n = p$. Thanks to the quantile formula \eqref{eq:W2_1d} for the expression of the Wasserstein metric when $d=1$, it follows that the  empirical Wasserstein barycenter  $\hat{\bmu}_{n,p}$ minimizing \eqref{eq:pbemp2Bar} is the unique discrete measure
\begin{equation}
\hat{\bmu}_{n,p} =  \frac{1}{p}\sum_{j=1}^p\delta_{\bar{\bX}_{j}^{\ast}}, \label{eq:expres1d}
\end{equation}
where   $\bX_{i,1}^{\ast} \leq \bX_{i,2}^{\ast} \leq \ldots \leq \bX_{i,p}^{\ast} $ are the order statistics the $i$-th sample, and $
\bar{\bX}_{j}^{\ast} = \frac{1}{n} \sum_{i=1}^{n} X_{i,j}^{\ast}$ for all $1 \leq j \leq p$. Moreover, as shown in  \cite{BGKL18},  the population Wasserstein barycenter $\mu_{\P}$ is the unique a.c.\ measure with quantile function
$$
F^{-}_{0}(\alpha) = \E \left( \bF^{-} (\alpha) \right) \mbox{ where $\bF^{-}$ is the quantile function of the random measure $\bnu$ with distribution $\P$.} 
$$
In this setting, the following upper bound on the expected quadratic risk $\E \left( W_{2}^{2} \left(\hat{\bmu}_{n,p}, \mu_{\P}  \right) \right)$ has been obtained in \cite{BGKL18}.  

\begin{thm} 
Let $Y_1,\ldots,Y_p$ be iid variables sampled from $\mu_{\P}$, and
$
\bmu_{p} = \frac{1}{p}\sum_{k=1}^{p} \delta_{Y_{k}}
$
Then, the estimator  $\hat{\bmu}_{n,p}$ satisfies
\begin{equation}
\E \left( W^{2}_{2}(\hat{\bmu}_{n,p}, \mu_{\P}) \right) \leq \frac{1}{n}  \int_{0}^{1}  \var\left(  \bF^{-}  (\alpha) \right)d\alpha +  \E \left( W^{2}_{2}(\bmu_{p}, \mu_{\P}) \right). \label{eq:upperbound1d}
\end{equation}
\end{thm}
Hence, the rate of convergence of $\hat{\bmu}_{n,p}$ is decomposed into the sum of two terms. The first one goes to zero at the parametric rate $\mathcal{O}(\frac{1}{n})$, while the second one depends on the rate of convergence of the empirical measure $\bmu_{p}$ towards it population counterpart in expected squared Wasserstein distance. For $d=1$, the question of deriving upper and lower bound for $\E \left( W^{2}_{2}(\bmu_{p}, \mu_{\P}) \right)$ has been thoroughly studied in \cite{W1}. In arbitrary dimension $d$, there is a rich literature in probability on studying the convergence rate of an empirical measure  towards it population version in Wasserstein distance. We refer to \cite{W1,fournier:hal-00915365,BachWeed,dedecker2019,BerthetWeed} for extensive references and the most recent contributions on this topic. For example, using results from  \cite{W1},  one has that
$$
\E \left( W_2^2(\bmu_{p}, \mu_{\P}) \right) \leq \frac{2}{p+1} J_{2}(\mu_{\P}), \mbox{ where }  J_{2}(\mu_{\P}) = \int_{\Omega} \frac{F_{0}(x)(1-F_{0}(x))}{f_{0}(x)} dx,
$$
where $F_0$ is the cumulative distribution function of $\mu_{\P}$ and $f_0$ denotes its density. Hence, if $J_{2}(\mu_{\P}) < + \infty$, it follows that $\hat{\bmu}_{n,p}$ converges to zero at the parametric rate $\mathcal{O}(\frac{1}{n} + \frac{1}{p})$. A detailed discussion on the derivation of upper bounds on the rate of converge of $\hat{\bmu}_{n,p}$  is proposed in \cite{BGKL18} with some results on minimax optimality of such bounds.



\subsubsection{Open problems in dimension $d \geq 2$}

Beyond the one-dimensional case $d=1$, it appears that deriving the rate of convergence for  empirical Wasserstein barycenters, e.g.\ for the expected quadratic risk $\E \left( W_{2}^{2} \left(\hat{\mu}_{\P_n}, \mu_{\P}  \right) \right)$,  is a much more involved problem. A few works \cite{spok2019,gouic2018} exist  on the convergence rate of  $\mu_{\P_{n}}$ to $\mu_{\P}$. However, the more general question of deriving the convergence rate of a Wasserstein barycenter $\hat{\bmu}_{n,p}$   \eqref{eq:pbemp2Bar} obtained from samples of random measure  is open for $d \geq 2$ (to the best of our knowledge).

\subsection{Regularization of Wasserstein barycenters}

An empirical Wasserstein barycenter $\hat{\bmu}_{n,p}$ of the discrete measures $\bnu_{p_{1}},\ldots,\bnu_{p_{n}}$, defined by \eqref{eq:discrete}, is generally irregular (and even not unique). Irregularity has to be understood in the sense that  $\hat{\bmu}_{n,p}$ is, in most situations, a discrete probability measure as shown by its expression \eqref{eq:expres1d} when $d=1$. Therefore, as a discrete measure, it poorly represents the smoothness of the population Wasserstein barycenter of a random  a.c.\ measure $\bnu$  with distribution $\P$, as $\mu_{\P}$ is absolutely continuous in this setting as discussed previously. Hence, we now explain how regularizing the  Wasserstein barycenter $\hat{\bmu}_{n,p}$ in order to construct a.c.\ and consistent estimators of an a.c.\ population barycenter in the asymptotic setting where both $n$ and $p = \min_{1 \leq i \leq n} p_{i}$ tend to infinity. To this end, one may consider two ways of regularizing an empirical barycenter.

\subsubsection{Penalized Wasserstein barycenters}

A first possibility is to follow the approach proposed in \cite{BCP17}. Thus, we now introduce the penalized Wasserstein barycenter associated to the distribution $\P$ that is defined as a solution of the minimization problem
$$ \min_{\mu\in\PP_2(\Omega)}\ \int_{\PP_2(\Omega)}W_2^2(\mu,\nu)d\P(\nu) +\gamma E(\mu)$$
where $\gamma > 0$ is a penalization parameter and the penalty function writes
\begin{equation}
\label{ex_penalty}
E(\mu) = \left\{\begin{array}{ll}
\Vert f\Vert_{H^{k}(\Omega)}^{2}, & \mbox{if}\ f =\frac{d\mu}{dx} \ \mbox{and} \ f\geq \alpha, \\
+\infty & \mbox{otherwise.}
\end{array}\right.
\end{equation}
where $\Vert\cdot\Vert_{H^k(\Omega)}$ denotes the Sobolev norm associated to the $\mathbb{L}^2(\Omega)$ space,  $\alpha>0$ is arbitrarily small and $k>d-1$. The choice \eqref{ex_penalty} for the penalty function $E$  is motivated by the discussion in Section 5 of \cite{BCP17}. Its main advantages are to impose that the penalized Wasserstein barycenter  is absolutely continuous with a smooth density.   Other examples of penalty functions are discussed in \cite{BCP17} including the class of relative G-functional described in Section 9.4 in \cite{ambrosio2008gradient}.

\begin{defi}\label{def:pen_bar}
For $\gamma>0$, we consider the Wasserstein barycenters:
\begin{align}
\mu^{\gamma}_{\P} &= \uargmin{\mu\in\PP_2(\Omega)} \ \int_{\PP_2(\Omega)}W_2^2(\mu,\nu)d\P(\nu)+\gamma E(\mu) \label{def:population_Regbar}\\
\hat{\bmu}^{\gamma}_{n,p} &= \uargmin{\mu\in\PP_2(\Omega)} \  \frac{1}{n}\sum_{i=1}^nW_2^2(\mu,\bnu_{p_{i}})+\gamma E(\mu) \label{def:simulated_Regbar}
\end{align}
called respectively penalized population barycenter \eqref{def:population_Regbar} and penalized empirical barycenter \eqref{def:simulated_Regbar}.
\end{defi}
Existence and uniqueness of these barycenters are discussed in \cite{BCP17} for a large class of penalty functions enforcing them to be a.c. One can thus let  $\hat{\bfun}^{\gamma}_{n,p}$ (resp.\ $f^{\gamma}_{\P}$) be  the densities  associated to $\hat{\bmu}^{\gamma}_{n,p}$ (resp.\  $\mu^{\gamma}_{\P}$). The choice \eqref{ex_penalty} for the penalty function $E$ imposes that the densities of the penalized Wasserstein barycenters are bounded from below by a small constant $\alpha > 0$ over their support that is equal to $\Omega$. In this setting, the following result (see Section 5 in \cite{BCP17}) gives the convergence rate in  expected squared  $\L_2(\Omega)$-distance between these densities for a fixed value of $\gamma > 0$.

\begin{thm}[Section 5 in \cite{BCP17}]\label{th:rate_Regbar}
Assume that $\Omega\subset\R^d$ is compact. Let $\hat{\bfun}_{n,p}^{\gamma}$ and $f_{\P}^{\gamma}$ denote the density functions of $\hat{\bmu}_{n,p}^{\gamma}$ and $\mu_{\P}^{\gamma}$, induced by the choice \eqref{ex_penalty} of the penalty function $E$. Then, provided that $d<4$, there exists a constant $c>0$ depending only on $\Omega$ such that
\begin{equation}\label{eq:rate_Regbar}
\E\left(\Vert \hat{\bfun}_{n,p}^{\gamma}-f_{\P}^{\gamma}\Vert^2_{\mathbb{L}_2(\Omega)}\right)\leq  c \left(\frac{1}{\gamma p^{1/4}}+\frac{1}{\gamma n^{1/2}}\right)
\end{equation}
where $p = \min_{1 \leq i \leq n} p_{i}$.
\end{thm}

Letting $f_{\P}^{0}$ be the density of the population Wasserstein barycenter $\mu_{\P}$,  the convergence of the approximation error term $\Vert f_{\P}^{\gamma}-f_{\P}^{0}\Vert^2_{\mathbb{L}_2(\Omega)}$ is also discussed in \cite{BCP17}, where it is shown that it converges to zero as $\gamma \to 0$. However, no rate of convergence (as a function of $\gamma$) has been derived in  \cite{BCP17}  for this approximation term, and deriving such a result remains an open problem for regularized  Wasserstein barycenters. Therefore, only the consistency of $\hat{\bfun}_{n,p}^{\gamma}$  has been derived in \cite{BCP17}  in the sense that the quadratic risk $\E\left(\Vert \hat{\bfun}_{n,p}^{\gamma}-f_{\P}^{0}\Vert^2_{\mathbb{L}_2(\Omega)}\right) \to 0$ under the assumption that $\gamma = \gamma_{n,p}$ decays to zero as $\min(n,p) \to + \infty$ at a rate which guarantees that the right hand size of inequality \eqref{eq:rate_Regbar} converges to zero. In Theorem \ref{th:rate_Regbar}, the condition $d<4$ may be relaxed at the price of obtaining a slower rate of convergence as $p \to + \infty$ (see Section 5 in \cite{BCP17}).

\begin{rem}
Consistency of regularized empirical Wasserstein barycenter is discussed in \cite{BCP17} for the Bregman divergence associated to the penalization function $E$. The specific choice \eqref{ex_penalty} for $E$ finally yields to the use of the $\mathbb{L}_2(\Omega)$ norm to compare the densities of  $\hat{\bfun}_{n,p}^{\gamma}$, $f_{\P}^{\gamma}$ and $f_{\P}^{0}$. The question of deriving consistency results for $\hat{\bmu}^{\gamma}_{n,p}$ using the Wasserstein metric (that is in the space of probability measures and not densities) has not been considered in \cite{BCP17}
\end{rem}

\subsubsection{Sinkhorn barycenters}

Another way to regularize an empirical Wasserstein barycenter is to use the Sinkhorn divergence \eqref{def:primal_Sinkhorn} to compare probability measures. As this divergence is defined for discrete measures supported on the same points, one first needs to choose a fixed grid $\Omega_{N}=\{x_1,\ldots,x_N\}$ made of $N$ points $x_k\in\mathbb{R}^d$ (bin locations that are typically equally spaced). Then, one performs a binning of the data \eqref{eq:discrete} on this grid, leading to a dataset of discrete measures (with supports included in $\Omega_{N}$) that we denote
\begin{equation}
\label{eq:data_bin}
\tilde{\bq}_i^{p_i}=\frac{1}{p_i}\sum_{j=1}^{p_i}\delta_{\tilde{X}_{i,j}}, \mbox{ where } \tilde{X}_{i,j} = \uargmin{x \in \Omega_{N}} \vert x - X_{i,j} \vert,
\end{equation}
for $1 \leq i \leq n$. Binning (i.e. choosing the grid $\Omega_{N}$) surely incorporates some sort of additional regularization, and a discussion on the influence of the grid size $N$ on the smoothness of the barycenter can be found in \cite{bigot2018datadriven}.   The size $N$ of the grid is also guided by numerical issues on the computational cost of the algorithms used to approximate a Sinkhorn divergence.

Then, another possibility to introduce regularization in the computation of  Wasserstein barycenters \cite{bigot2018datadriven,cuturi2013fast,CuturiPeyre} is to consider the following estimator
\begin{equation}
\hat{\br}^{\varepsilon}_{n,p} =\uargmin{r\in\Sigma_N}\frac{1}{n}\sum_{i=1}^n W_{2,\varepsilon}^2(r,\tilde{\bq}_i^{p_i})  , \label{eq:defr}
\end{equation}
that can be interpreted as a  Fr\'echet mean with respect to a Sinkhorn divergence. We call $\hat{\br}^{\varepsilon}_{n,p}$ a Sinkhorn barycenter. From the point of view of statistics, it can be interpreted as a M-estimator defined over the compact set $\Sigma_N$. Two key properties of the Sinkhorn divergence have been used in \cite{bigot2018datadriven} to study the convergence properties of such barycenters:
\begin{description}
\item[(i)]  for any $q\in\Sigma_N$, the function $r\mapsto W_{2,\varepsilon}^2(r,q)$ is  $\varepsilon$- strongly convex for the Euclidean 2-norm 
\item[(ii)]  let $q \in \Sigma_N$ and $0 < \rho < 1/N$ an arbitrarily small constant. Then, one has that $r  \mapsto W_{2,\varepsilon}^2(r,q)$ is $L_{\rho,\varepsilon}$-Lipschitz on
$$
\Sigma_N^{\rho} = \left\{ r \in \Sigma_N \; : \; \min_{1 \leq \ell \leq N} r_{\ell} \geq \rho \right\},
$$
with
$$
L_{\rho,\varepsilon}= \left( \sum_{1 \leq m \leq N} \left(2 \varepsilon \log(N) + \sup_{1 \leq \ell, k \leq N} |C_{m\ell} - C_{k \ell}|    - 2 \varepsilon \log(\rho) \right)^2  \right)^{1/2}
$$
\end{description}

For a discussion on the behavior of the Lipschitz constant $L_{\rho,\varepsilon}$ as $\rho$ and $\varepsilon$ tend to zero, we refer to Section 3.1 in  \cite{bigot2018datadriven}. For example, if $\rho = \varepsilon^{\kappa}$ for some $\kappa > 0$, then $$\lim_{\varepsilon \to 0} L_{\rho,\varepsilon} = \left( \sum_{1 \leq m \leq N} \left( \sup_{1 \leq \ell, k \leq N} |C_{m\ell} - C_{k \ell}|  \right)^2  \right)^{1/2}.$$ To guarantee the Lipschitz continuity of the mapping $r\mapsto W_{2,\varepsilon}^2(r,q)$, the analysis of the statistical properties of Sinkhorn barycenters in \cite{bigot2018datadriven} is restricted to discrete measures belonging the convex set $\Sigma_N^{\rho}$ (thus having non-vanishing entries). Therefore, one has to slightly modify the data at hand so that the observations belong to $\Sigma_N^{\rho}$. This yields to the following definitions of sampling model, and empirical/population Sinkhorn barycenters in  $\Sigma_N^{\rho}$.

\begin{defi}\label{def:Sinkhorn_bar}
Let $0 < \rho < 1/N$, and $\P$ be a probability distribution on $\Sigma_N^{\rho}$. Let $\bq_1,\ldots,\bq_n\in\Sigma_N^{\rho}$ be an iid sample drawn from the distribution $\P$. Assume that $(\tilde{X}_{i,j})_{1\leq j\leq p_i}$ are iid random variables sampled from $\bq_i$,  for each $1 \leq i \leq n$ and consider the  discrete measures
$
\tilde{\bq}_i^{p_i} = \frac{1}{p_i}\sum_{j=1}^{p_i}\delta_{\tilde{X}_{i,j}}
$
and
$
\hat{\bq}_i^{p_i}  = (1 - \rho N)\tilde{\bq}_i^{p_i} +  \rho \1_{N},
$
where $\1_{N}$ is the vector of $\R^N$ with all entries equal to one. Then, we define
\begin{align}
r^{\varepsilon} &=\uargmin{r\in\Sigma_N^{\rho}}\E_{\bq\sim\P}[W_{2,\varepsilon}^2(r,\bq)] & \qquad \mbox{the population Sinkhorn  barycenter}\label{def:pop_emp_Sinkhorn}\\
\hat{\br}^{\varepsilon}_{n,p} &=\uargmin{r\in\Sigma_N^{\rho}}\frac{1}{n}\sum_{i=1}^nW_{2,\varepsilon}^2(r,\hat{\bq}_i^{p_i}) & \qquad \mbox{the empirical Sinkhorn  barycenter}\label{def:simulated_Sinkhorn}
\end{align}
\end{defi}

\begin{thm}[Section 3 in \cite{bigot2018datadriven}]\label{th:rate_Sinkbar}
Let $\varepsilon >0$ and $p=\min_{1 \leq i \leq n} p_i$. Then, one has that
\begin{equation}
\E(\| r^{\varepsilon}-\hat{r}^{\varepsilon}_{n,p} \|^2_{\R^N}) \leq  \frac{32 L_{\rho,\varepsilon}^2}{\varepsilon^2n}  +  \frac{2 L_{\rho,\varepsilon}}{\varepsilon} \left(   \sqrt{\frac{N}{ p}} + 2 \rho( N + \sqrt{N} ) \right). \label{eq:rateSinkbar}
\end{equation}
\end{thm}

The upper bound \eqref{eq:rateSinkbar} allows to shed some light on how the expected squared distance between the empirical Sinkhorn  barycenter and its population version is influenced by the number $n$ of observed measures, the minimal number $p$ of samples per measure and the size $N$ of the grid. However, this upper bound does not converge to zero as $\min(n,p) \to + \infty$ if $\rho$ remains fixed. The main advantage of this upper bound is to allow the derivation of a data-driven strategy to select $\varepsilon$ that is presented in the following subsection below. Finally, one should remark that the behavior of the population Sinkhorn  barycenter $r^{\varepsilon}$  as $\varepsilon \to 0$ has not been considered in \cite{bigot2018datadriven}.

\begin{rem}
A third way to regularize an empirical Wasserstein barycenter (that might look more natural for statisticians) is to perform a preliminary smoothing step of the discrete measures $\bnu_{p_{1}},\ldots,\bnu_{p_{n}}$. For example, one may first smooth the data using standard kernel density estimation
$
\hat{\bfun}^{h_i}_{i}(x) = \frac{1}{p_{i} h_{i}^d} \sum_{j=1}^{p_{i}} K\left( \frac{x-\bX_{i,j}}{h_{i}}\right)
$
and then define a smoothed Wasserstein barycenter as
$$
\hat{\bmu}^{h}_{n,p} = \argmin_{\mu \in \PP_2(\Omega)} \frac{1}{n} \sum_{i=1}^{n} W_{2}^{2}(\hat{\bnu}_{i}^{h_{i}},\mu),
$$
where $d \hat{\bnu}_{i}^{h_{i}}(x) = \hat{\bfun}^{h_i}_{i}(x) dx $, and $h = \max_{1 \leq i \leq} h_i$. The consistency of this approach (as $h \to 0$) has been studied in details by \cite{Pana15,Pana18} using multiple point processes for the sampling procedure, and it is also discussed in \cite{BGKL18}.
\end{rem}

\subsection{Data-driven regularization in computational optimal transport}

In this paper, we report results on numerical experiments from \cite{bigot2018datadriven} that are based on the approach proposed in \cite{CuturiPeyre} to compute Sinkhorn barycenters on a grid $\Omega_{N}$ of equi-spaced points. We do not report results on the computational aspects  of penalized Wasserstein barycenters, and we refer to  \cite{bigot2018datadriven} for illustrations on their numerical performances. 

\begin{figure}
\centering
\subfigure[]{\includegraphics[width=0.95 \textwidth,height=0.45\textwidth]{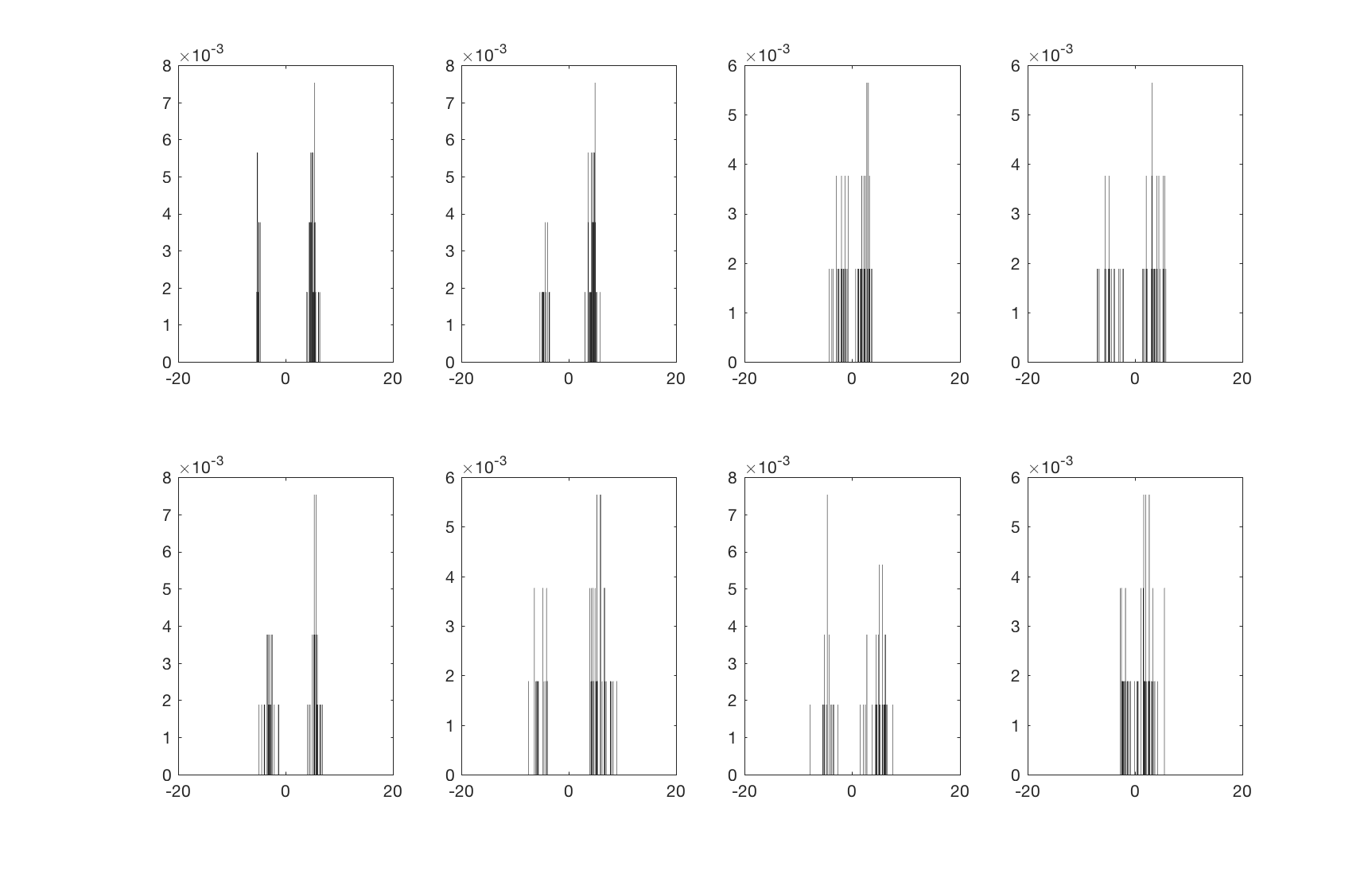}}
\subfigure[]{\includegraphics[width=0.45 \textwidth,height=0.35\textwidth]{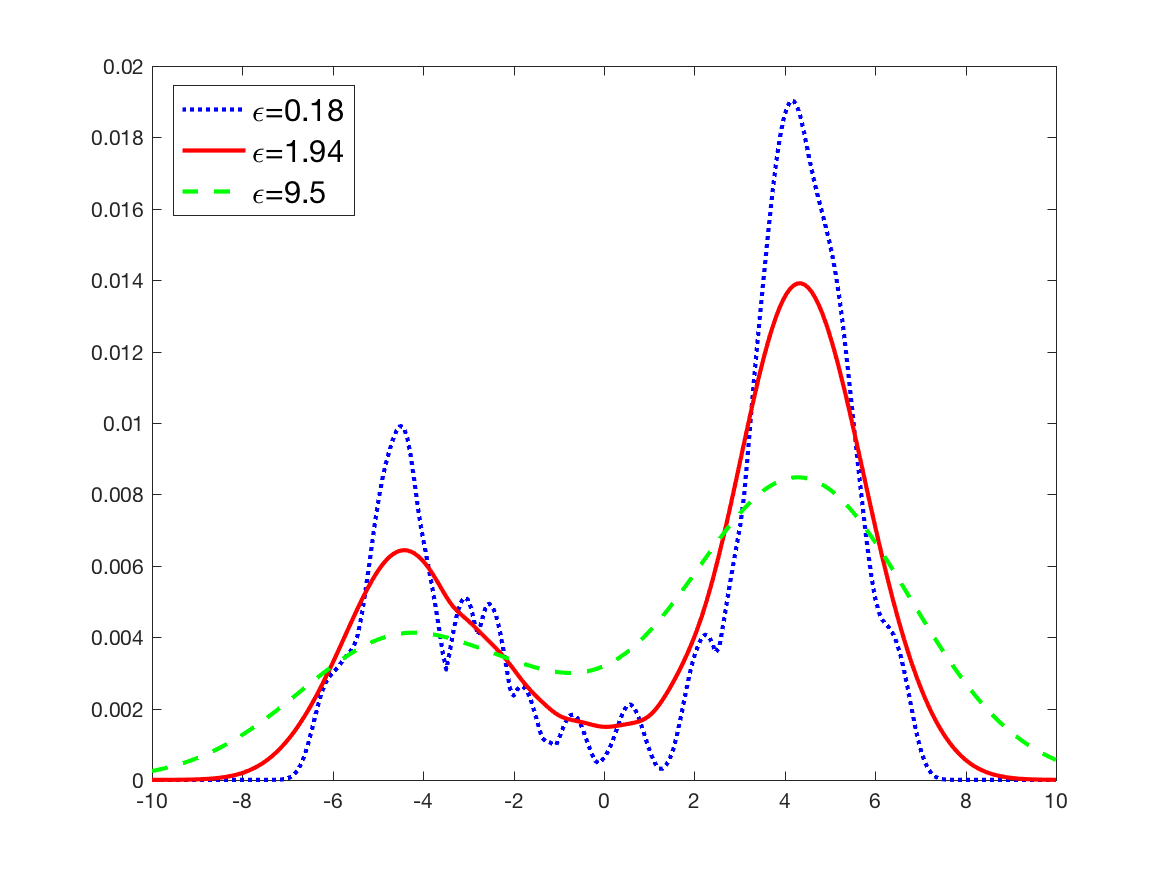}}
\caption{(a) A subset of $8$ histograms (out of $n=15$) obtained with random variables sampled from  one-dimensional Gaussian mixtures distributions $\bnu_i$ (with random means and variances). The population Wasserstein of $\bnu_1,\ldots,\bnu_n$ is expected to be a mixture of two Gaussian distributions. Histograms are constructed by binning the data  $(X_{i,j})_{1\leq i\leq n\ ; 1\leq p}$ on a grid $\Omega_N$ of size $N=2^8$ with $p_1 =\ldots = p_n = p = 50$. (b) Three Sinkhorn barycenters $\hat{\br}_{n,p}^{\hat{\varepsilon}}$ associated to the parameters $\varepsilon=0.18, 1.94, 9.5$. }
\label{fig:ex_gaussian1D}
\end{figure}

As any nonparametric method, the use of regularized Wasserstein barycenters involves the delicate calibration of a regularization parameter. Choosing the value of $\varepsilon$ for the Sinkhorn barycenter $\hat{r}^{\varepsilon}_{n,p}$ might be guided by the classical bias and variance tradeoff principle as illustrated by the synthetic example  displayed in Figure \ref{fig:ex_gaussian1D}. For small values of $\varepsilon$ the Sinkhorn barycenter shows many oscillations while it is much smoother for larger values of $\varepsilon$. In some sense, the regularization parameters $\varepsilon$ can be interpreted as the usual bandwidth parameter in kernel density estimation and its  choice greatly influences the shape of $\hat{r}^{\varepsilon}_{n,p}$.

\begin{figure}
\centering
\subfigure[]{\includegraphics[width=0.45 \textwidth,height=0.45\textwidth]{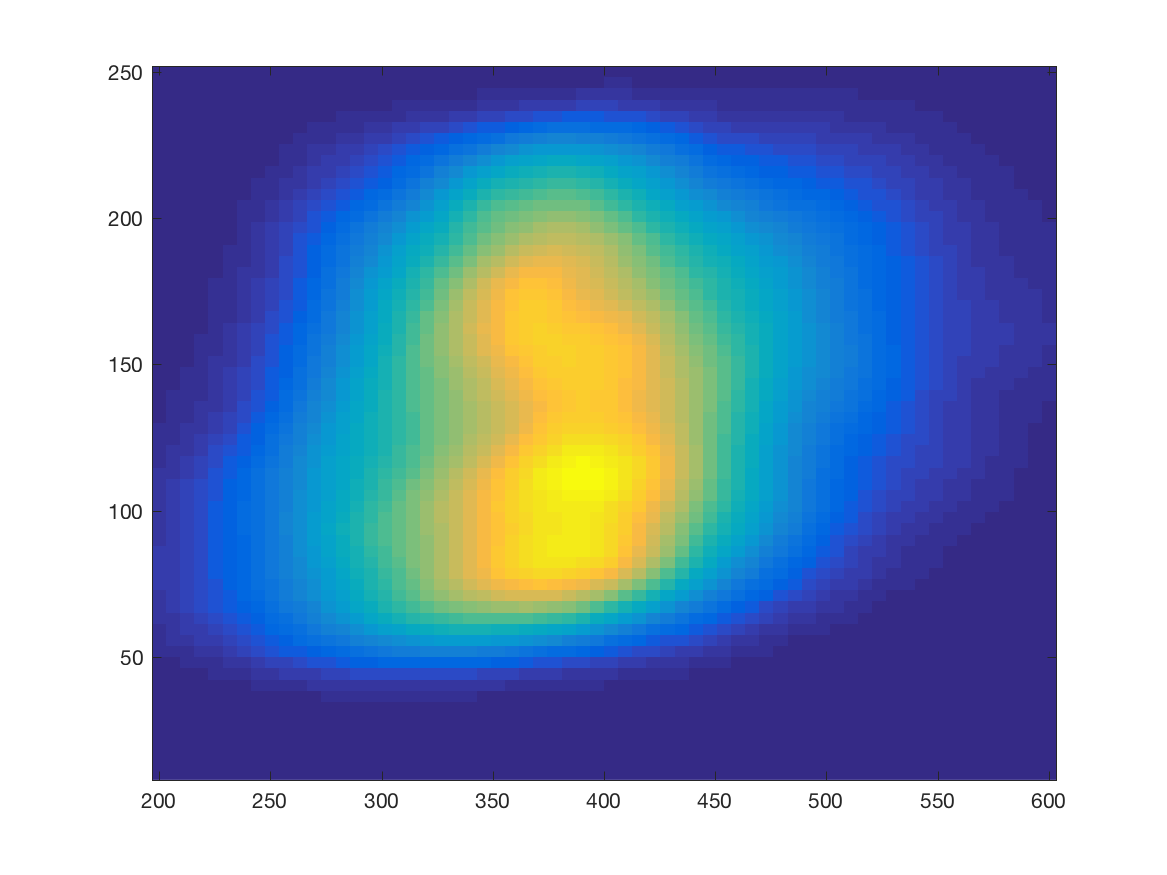}}
\subfigure[]{\includegraphics[width=0.45 \textwidth,height=0.45\textwidth]{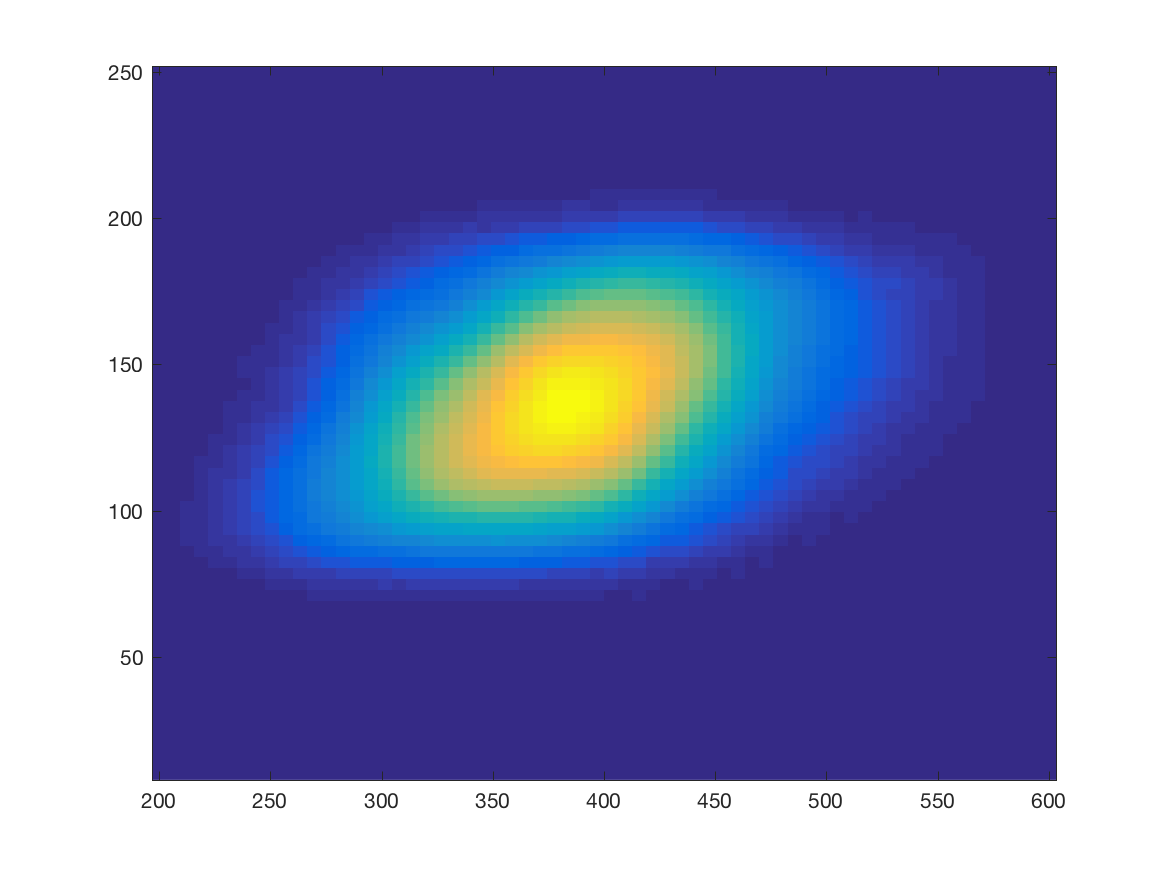}}
\caption{Two dimensional flow cytometry dataset. (a) Euclidean mean  of the data (after smoothing but without a preliminary step to correct misalignment effects in the raw data), (b)  Sinkhorn barycenter $\hat{\br}_{n,p}^{\hat{\varepsilon}}$ associated to the parameter $\hat{\varepsilon}=2.6$.}
\label{fig:mean_cytometry2D}
\end{figure}

Let us now briefly discuss the question of choosing  $\varepsilon$ in an automatic way.
In \cite{bigot2018datadriven}, the following (un-regularized) population Wasserstein is considered
\begin{equation}\label{def:r0}
r^{0} \in \uargmin{r\in\Sigma_N^{\rho}}\E_{\bq\sim\P}[W_{2,0}^2(r,\bq)] \quad \mbox{with} \quad
W_{2,0}^2(r,q):=\underset{U\in U(r,q)}{\min}\  \langle C, U\rangle.
\end{equation}
The term $\E(\vert r^{\varepsilon}-\hat{r}^{\varepsilon}_{n,p}\vert^2)$ can then be interpreted as ``variance term'', while  $\vert r^{\varepsilon} - r^{0}\vert$ is referred to as a bias term. Ideally, one would like to select a value of $\varepsilon$ which gives  a good tradeoff between these two terms whose values change in opposite directions as $\varepsilon$ varies. However, the decay of the bias term (as $\varepsilon \to 0$) is typically unknown and only an estimation of the variance term is available.  In \cite{bigot2018datadriven}, it is thus proposed to use the upper bound  \eqref{eq:rateSinkbar} to automatically calibrate the parameter $\varepsilon>0$ by following Goldenshluger-Lepski (GL) principle \cite{goldenshluger2008universal} as formulated in  \cite{lacour2016minimal}  for standard kernel density estimation. In \cite{bigot2018datadriven}, an adaptation of the GL's principle is used to  to estimate a bias-variance trade-off functional that is minimized over a finite set of regularization parameters to obtain an automatic selection of $\varepsilon$ for Sinkhorn barycenters. The method consists in comparing estimators pairwise (for a given range of regularization parameters) with respect to a loss function, and we refer to \cite{bigot2018datadriven} for further details. As an illustration of the usefulness of the GL's principle in statistical optimal transport, we report in Figure \ref{fig:mean_cytometry2D} numerical experiments from \cite{bigot2018datadriven} on the computation of a Sinkhorn barycenter (with a data-driven choice of $\varepsilon$) in dimension $d=2$ for the flow cytometry dataset displayed in Figure \ref{fig:ex_cytometry2D}(b). This regularized barycenter clearly corrects some mis-alignment issues in these measurements. We also display in Figure \ref{fig:mean_cytometry2D}(a) the Euclidean mean of this dataset (after a preliminary kernel smoothing step). The support of this Euclidean mean is more  spread out due to the presence of a strong variance in spatial location in the data from one individual to another.

\section{Geometry of the Wasserstein space and geodesic PCA} \label{sec:GPCA}

Principal Component Analysis (PCA) is certainly the most widely used approach to reduce the dimension of datasets whose elements are high-dimensional vectors or matrices. When such elements can be modeled as probability measures, we explain how using the Wasserstein metric  to compute their modes of variation accordingly around their barycenter in the Wasserstein space. We briefly describe the setting where the data can be modeled as histograms supported on the real line, and we refer to extensions in the case of probability measures supported on $\R^d$ for $d \geq 2$. 

\subsection{Geometry of the Wasserstein space}

 The space $\PP_2(\Omega)$ endowed with the 2-Wasserstein distance is not a Hilbert space. As a consequence, standard PCA in a linear space  based on the diagonalisation of a covariance matrix cannot be applied directly to compute principal modes of variations with respect to the Wasserstein metric. Nevertheless, a meaningful notion of PCA \cite{NIPS2015_5680,bigot2017geodesic,bigot2018GPCA} can still be defined by relying on the pseudo-Riemannian structure of the Wasserstein space which was extensively studied in \cite{ambrosio2008gradient}. 
 
 This geometric structure allows to define a notion of geodesic PCA in the Wasserstein space which  shares similarities with analogs of PCA for data belonging to a Riemannian manifolds. Principal Geodesic Analysis (PGA) \cite{geodesicPCA,exactPGA}   is a geometric method for statistical analysis that extends the notion of classical PCA in Hilbert spaces \cite{MR650934,MR2168993} to curved Riemannian manifolds  by replacing Euclidean concepts of vector means, lines and orthogonality by the more general notions of Fr\'echet mean, geodesics in Riemannian manifolds and orthogonality in tangent spaces. For $d=1$, the geometric structure of the Wasserstein space is relatively easy to describe, and we refer to \cite{bigot2017geodesic,bigot2018GPCA} for  a detailed presentation. The setting $d \geq 2$ is more involved. Some elements on the geometry of the Wasserstein space can be found in Section 4.4 in \cite{Pana18}, and an extensive study can be found in \cite{ambrosio2008gradient}.

\subsection{Geodesic PCA in the one-dimensional case}

Following these principles of statistical analysis in manifolds, a framework for geodesic PCA (GPCA) of probability measures supported on a interval $\Omega \subset \R$ has been introduced in \cite{bigot2017geodesic}. GPCA is defined as the problem of estimating a principal geodesic subspace (of a given dimension) which maximizes the variance of the projection of the data to that subspace whose base point is the Wasserstein barycenter of the observations. In \cite{bigot2017geodesic}, a detailed study is proposed on the existence and consistency of GPCA. It is also shown that GPCA in the Wasserstein space is equivalent to map the data in the tangent space  at the empirical Wasserstein barycenter (assumed to be a.c.) and then to perform a PCA that is constrained to lie in a convex and closed subset of functions. Mapping the data to this tangent space is not difficult in the one-dimensional case as it amounts to computing a set of optimal maps between the data and their Wasserstein barycenter, for which a closed form is available using their quantile functions. In particular, optimal maps belong to the set of non-decreasing functions. This convex constrained PCA can also be understood from the quantile formula \eqref{eq:W2_1d} for the Wasserstein metric when $d=1$. Indeed, GPCA for probability measures supported on $\Omega \subset \R$ corresponds to a PCA of quantile functions that is constrained to lie in the convex subset of non-decreasing function in $\mathbb{L}_2([0,1])$. Nevertheless, such a convex PCA leads to the minimization of a non-convex and non-differentiable objective function. This is a delicate optimisation problem (even for $d=1$), and numerical methods to compute such a PCA have been proposed in \cite{bigot2018GPCA}. 

\begin{figure}
{\includegraphics[width=0.95 \textwidth,height=0.55\textwidth]{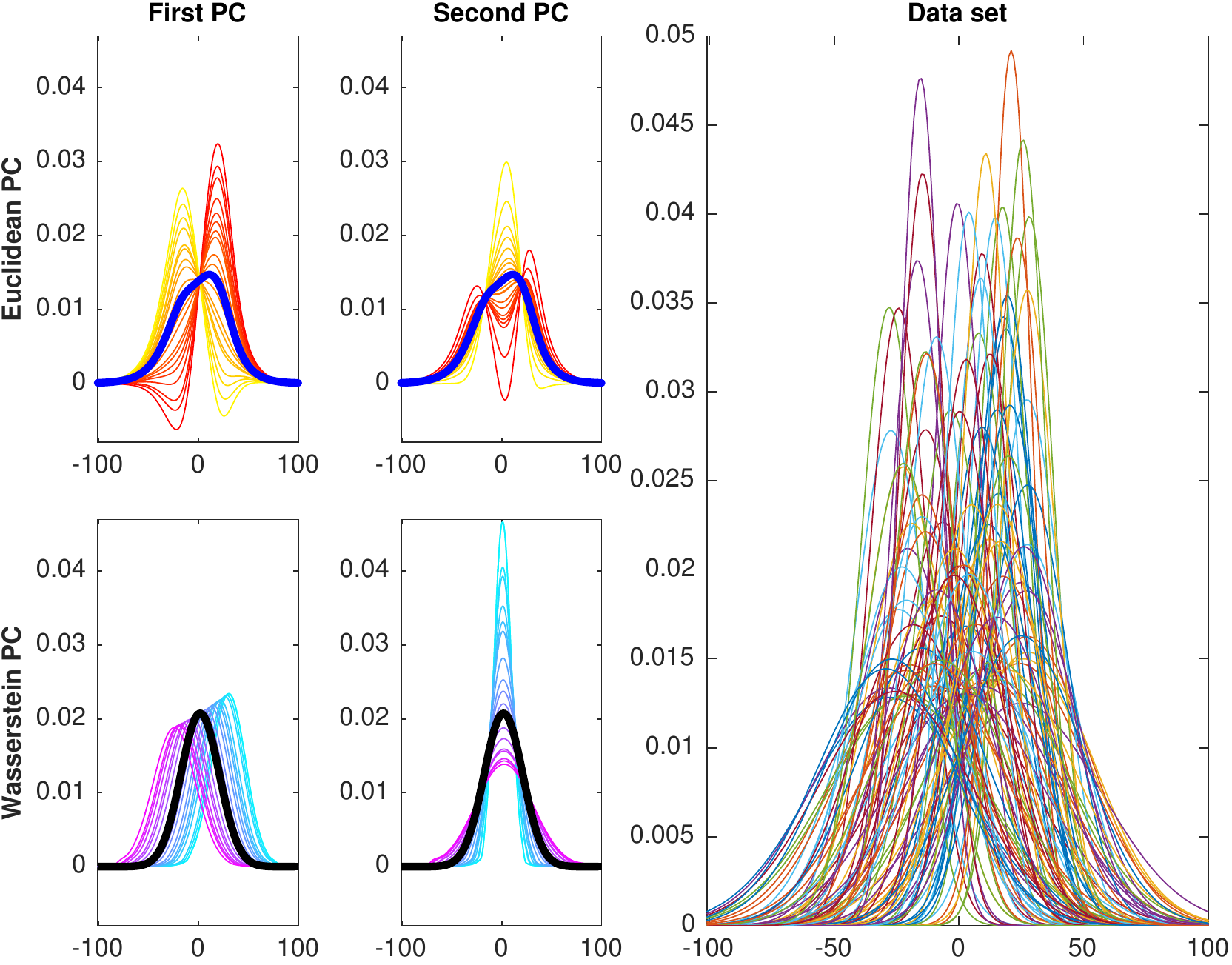}}
\caption{Synthetic data set: Gaussian probability densities with random mean and variance. {\bf Euclidean PC}: standard PCA with respect to the Euclidean metric between densities - first and second principal components (PC)  of the data around their Euclidean mean (blue curve). The yellow to red curves encode the progression of the first two main directions of variation in the Hilbert space of square integrable functions. {\bf Wasserstein PC}: geodesic PCA in the Wasserstein space - first and second principal components (PC)  of the data  around their Wasserstein barycenter (black curve). The curves from light blue to violet encode the progression of the pdf of the first two main geodesics directions of variation in $\PP_2(\Omega)$.} \label{fig:exPCA1D}
\end{figure}

\begin{figure}
\subfigure[]{{\includegraphics[width=0.32 \textwidth,height=0.32\textwidth]{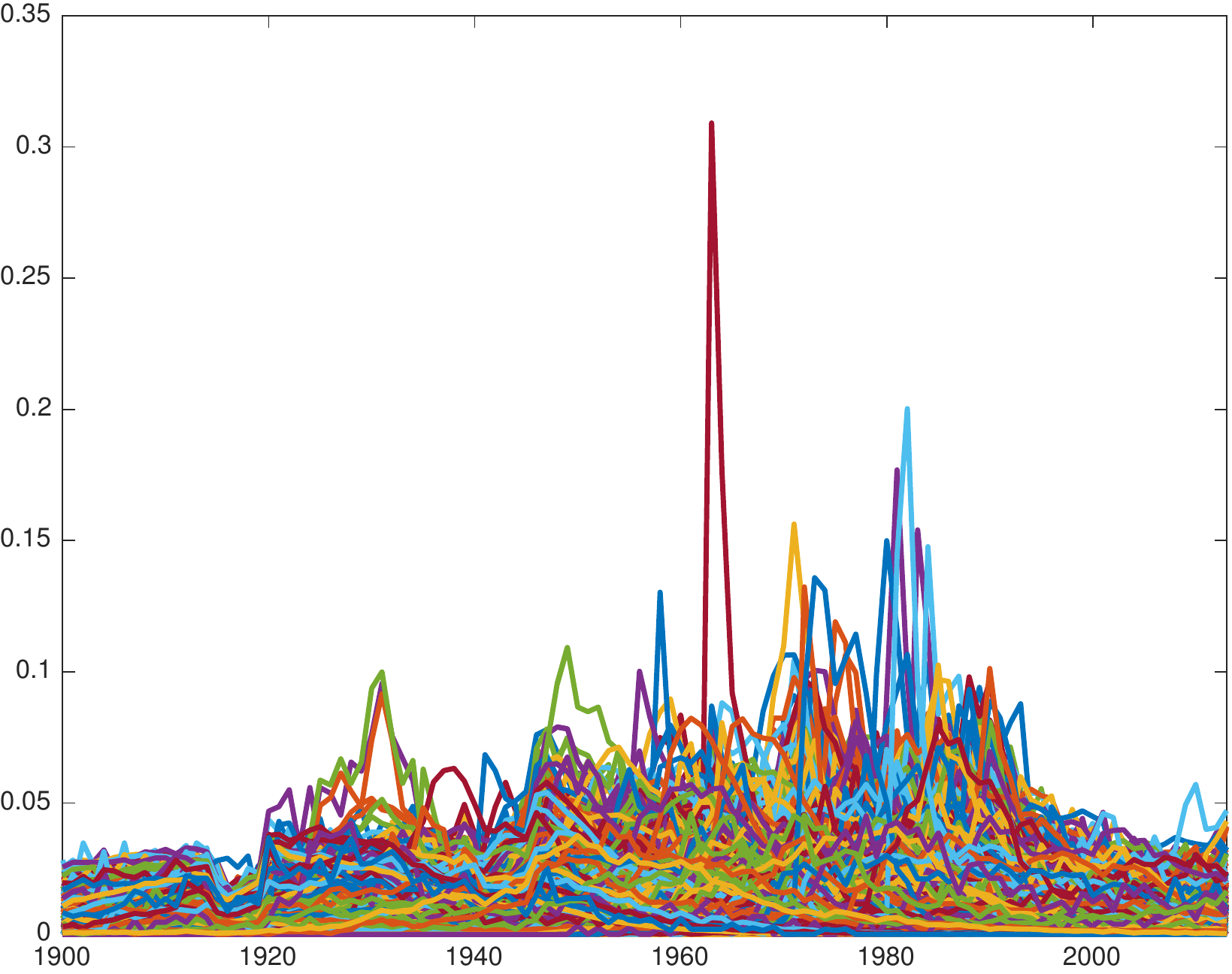}}}
\subfigure[]{{\includegraphics[width=0.32 \textwidth,height=0.32\textwidth]{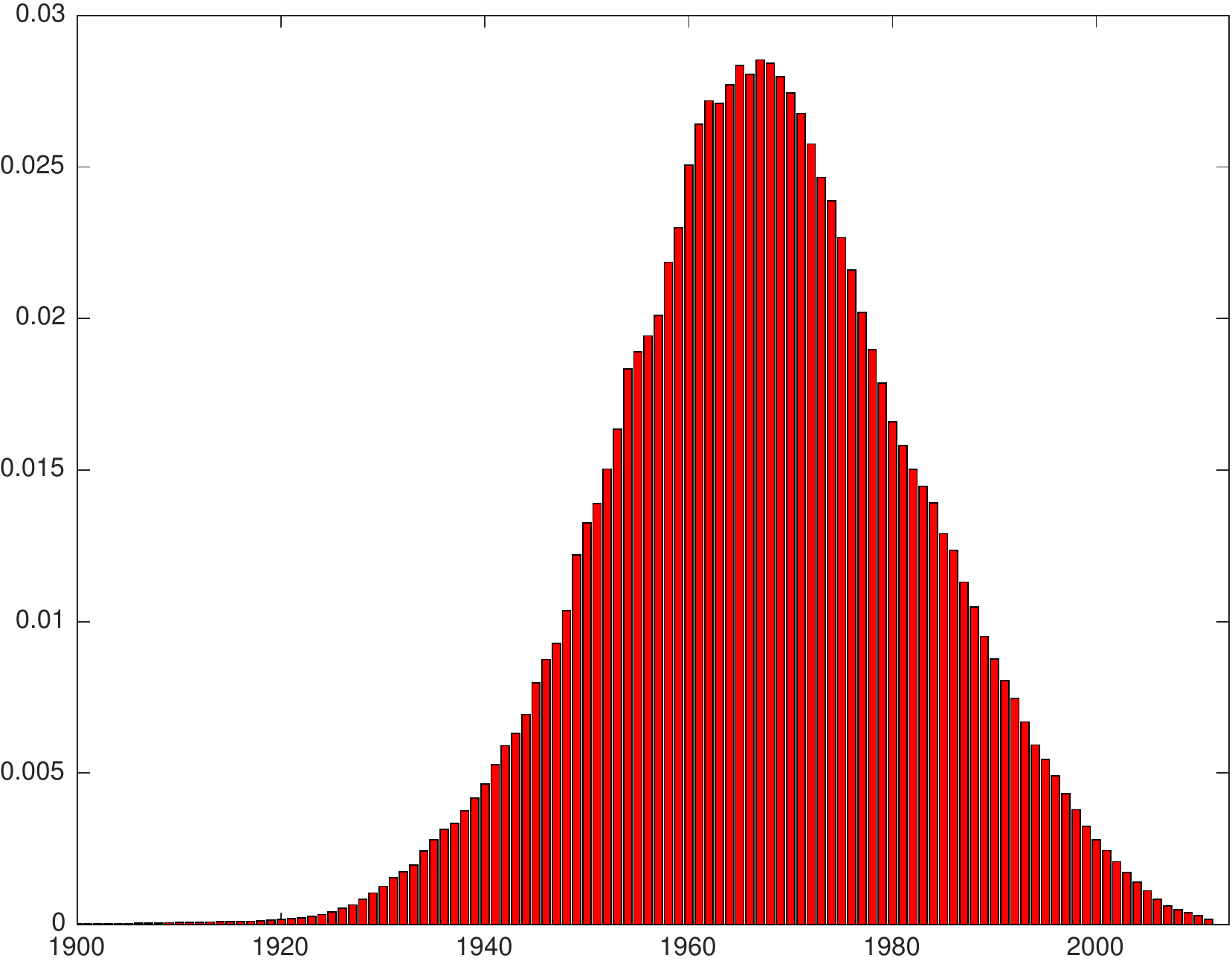} }}

\subfigure[]{{\includegraphics[width=0.32 \textwidth,height=0.32\textwidth]{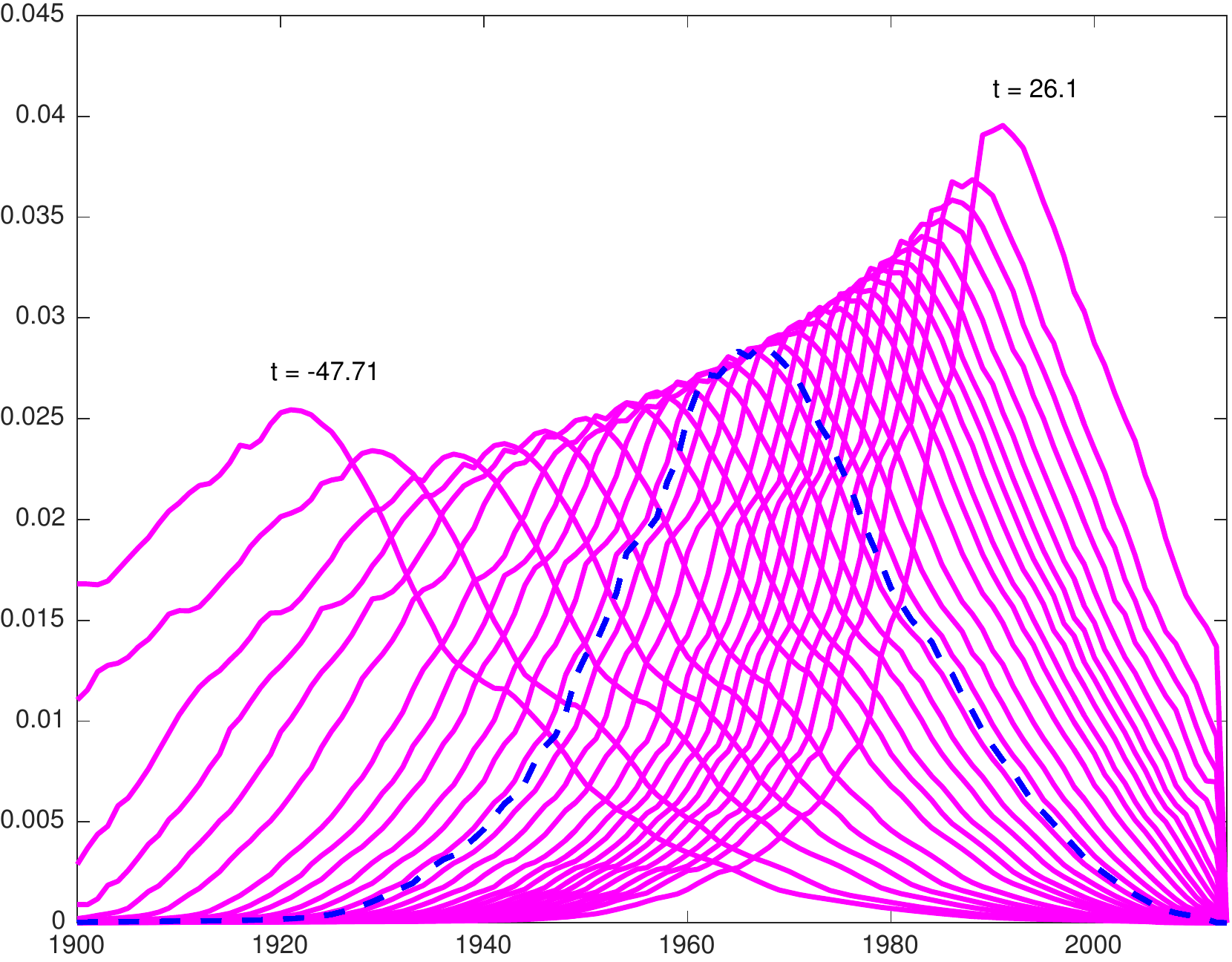}}}
\subfigure[]{{\includegraphics[width=0.32 \textwidth,height=0.32\textwidth]{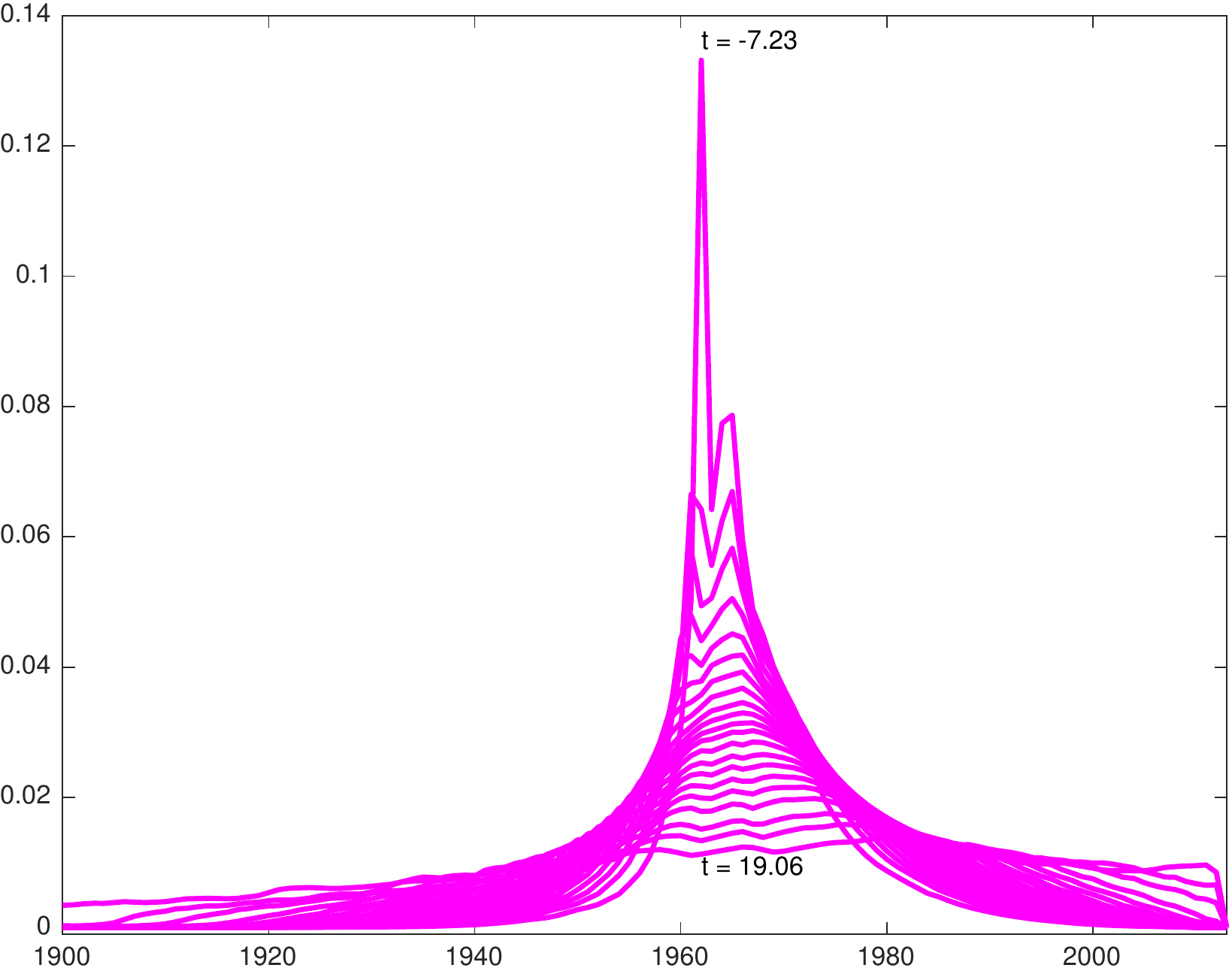}}}
\caption{Geodesic PCA of the data on children's first name at birth. (a) The whole dataset made of $n=1060$ histograms. (b) Density of the Wasserstein barycenter.  (c) First direction and (d) second direction  of geodesic variation in $\PP_2(\Omega)$ (at the level of pdf).} \label{fig:exPCAnames}
\end{figure}

To shed some lights on the benefits of GPCA for data analysis, we report numerical experiments on synthetic data in Figure \ref{fig:exPCA1D} and for the dataset on children's first name at birth in Figure \ref{fig:exPCAnames}. When the size or locations of significant bins in observed histograms significantly vary from one  histogram to another, standard PCA with respect to the Euclidean metric (in the space of densities)  fail in recovering geometric variations in location. To the contrary, GPCA in the Wasserstein space is able to capture meaningful sources of variations in a dataset by separating variability in location and scaling in the observed histograms along the first and second geodesic principal components.

\subsection{Geodesic PCA  in dimension $d \geq 2$}

The study of Geodesic PCA of histograms in the Wasserstein space is much more involved  in dimension $d \geq 2$, both on the theoretical and numerical aspects. Indeed, for $d \geq 2$, one could always define GPCA as the construction of a principal geodesic subspace maximizing the variance of the projection of the data to that subspace. However, a rigorous analysis of such a construction in dimension $d \geq 2$ has not been proposed so far. There also exist two main differences with the one-dimensional case. First,  GPCA cannot be equivalent to map the data in the tangent space  at their  barycenter (assumed to be a.c.) and then to perform a standard PCA in that linear space as  the Wasserstein space is a curved manifold for $d \geq 2$. Secondly, optimal maps sending a set of histograms to such a tangent space are more difficult to compute for $d \geq 2$. Moreover, it is well known  \cite{Brenier91,cuesta1989} that  such maps are equal to the gradient of convex functions. This implies that an adaptation of existing algorithms for $d=1$  leads to  optimisation problems with more difficult constraints to handle than a PCA constrained to lie in the convex set of non-decreasing functions.

\begin{figure}
{\includegraphics[width=0.75 \textwidth,height=0.55\textwidth]{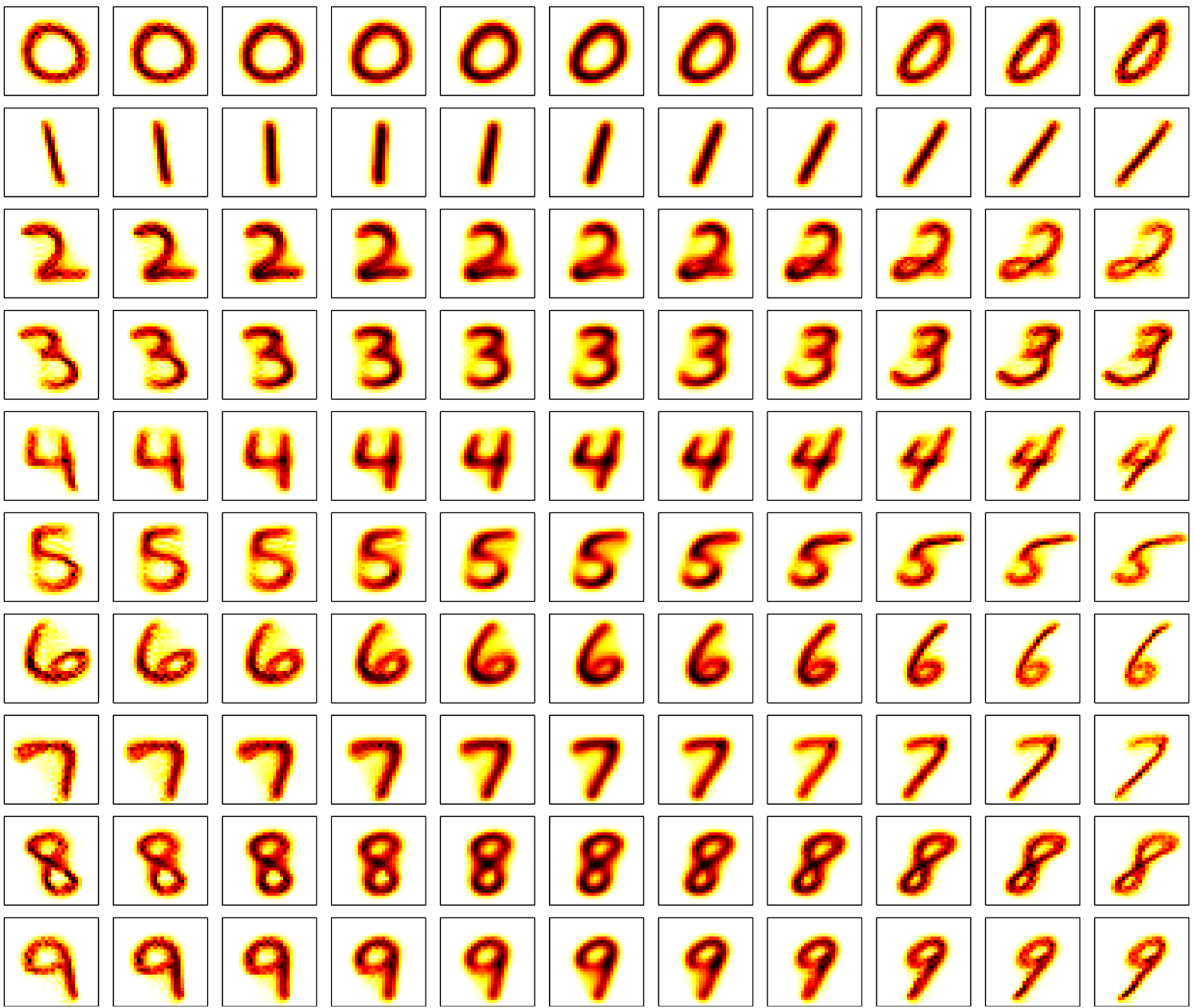}}
\caption{First principal geodesics for 1000 images of each digit from the MNIST database  \cite{lecun1998mnist}. The middle column represents the pdf of the Wasserstein barycenter of each digit (computed from 1000 images whose pixels are normalized to form a set of histograms). 
Each line allows to visualize the first geodesic of variation of each digit around their barycenter in the sense of  optimal transport. The cost is quadratic and it corresponds to the squared distance between pixel locations.} \label{fig:MNIST}
\end{figure}

Nevertheless, numerical approximation to the construction of GPCA in dimension $d=2$ have been proposed in \cite{NIPS2015_5680,bigot2018GPCA}. As an illustrative example, we report numerical experiments from \cite{bigot2018GPCA} on the analysis of data from the  MNIST database \cite{lecun1998mnist} which contains grayscale images of handwritten digits. All grayscale images in this dataset are of size $28 \times 28$ pixels. We normalize them so that the sum of pixel grayscale values sum to one. In this manner, the images are histograms supported on a 2D grid of size $28 \times 28$. The ground metric for the Wasserstein distance is then the squared Euclidean distance between the locations of bins in this 2D grid. In Figure \ref{fig:MNIST}, we display the first principal components computed from 1000 images of each digit using the algorithmic approach proposed in \cite{bigot2018GPCA}. It can be seen that  GPCA  captures very well the main source of variability of each digit such as changes in the shape of the '0' or the presence or absence of the lower loop of the '2'. 

\bibliographystyle{siam}
\bibliography{data_analysis_W2}

\end{document}